\documentclass[11pt]{article}

\usepackage[a4paper,
            includefoot,
            left=2cm, right=1.5cm, top=2cm, bottom=2cm, headsep=1cm, footskip=1cm]{geometry}
\usepackage[x11names]{xcolor}
\usepackage[unicode=true,pdfusetitle,
 bookmarks=true,bookmarksnumbered=false,bookmarksopen=true,bookmarksopenlevel=1,
 breaklinks=false,pdfborder={0 0 1},backref=false]
 {hyperref}
 \hypersetup{
    colorlinks = false,
    linkbordercolor  = green,
    citebordercolor = red,
 }

\usepackage{lipsum}
\usepackage{bm}
\usepackage{amsmath}
\usepackage{amssymb}
\usepackage{amsthm}
\usepackage{graphicx}
\usepackage{epstopdf}

\ifpdf
  \DeclareGraphicsExtensions{.eps,.pdf,.png,.jpg}
\else
  \DeclareGraphicsExtensions{.eps}
\fi

\usepackage{amsmath}
\numberwithin{equation}{section}

\usepackage{amsopn}
\DeclareMathOperator{\diag}{diag}
\DeclareMathOperator{\mes}{mes}

\usepackage{euscript}
\usepackage{relsize}
\usepackage{mathdots}
\usepackage{graphicx}
\usepackage{indentfirst}
\usepackage{empheq}
\usepackage{multirow}
\usepackage{thmtools}
\usepackage{thm-restate}
\usepackage{tikz}
\usepackage{algorithm}
\usepackage{algpseudocode}

\DeclareMathOperator*{\argmax}{arg\,max}
\DeclareMathOperator*{\argmin}{arg\,min}

\newcommand{\row}[1]{#1 ,\,\boldsymbol{:}\,}
\newcommand{\vecrow}[1]{#1}
\newcommand{\col}[1]{\,\boldsymbol{:}, #1}

\def\rank{\mathop{\mathrm{rank}}}

\def\colspace{\mathop{\mathrm{colspace}}}
\def\rowspace{\mathop{\mathrm{rowspace}}}
\def\Sigminus{\bfW}
\newcommand\winverse[2]{\left(#1 \right)^{\dag}_{#2}}
\newcommand\inverse[1]{#1^{\dag}}
\def\i0{\tau}
\def\Si0{\dot \bfs}
\def\Ai0{\dot \bfa}
\def\Proj{\mathbf{\Pi}}

\def\fullop{H_{\tau}}

\newcommand{\Input}{\hspace*{\algorithmicindent} \textbf{Input}:\ }
\newcommand\upangle[2]{\mathord{<\mspace{-9mu}\mathrel{)}\mspace{2mu}}(#1, #2)}

\newtheorem{proposition}{Proposition}[section]
\newtheorem{theorem}{Theorem}[section]
\newtheorem{remark}{Remark}[section]
\newtheorem{lemma}{Lemma}[section]
\usepackage{euscript}

\newcommand{\spC}{\mathbb{C}}

\newcommand{\spR}{\mathbb{R}}

\newcommand{\spT}{\mathbb{T}}

\newcommand{\tsA}{\mathsf{A}}
\newcommand{\tsB}{\mathsf{B}}
\newcommand{\tsC}{\mathsf{C}}
\newcommand{\tsD}{\mathsf{D}}

\newcommand{\tsP}{\mathsf{P}}

\newcommand{\tsR}{\mathsf{R}}
\newcommand{\tsS}{\mathsf{S}}

\newcommand{\tsX}{\mathsf{X}}
\newcommand{\tsY}{\mathsf{Y}}
\newcommand{\tsZ}{\mathsf{Z}}

\newcommand{\bfA}{\mathbf{A}}
\newcommand{\bfB}{\mathbf{B}}
\newcommand{\bfC}{\mathbf{C}}
\newcommand{\bfG}{\mathbf{G}}
\newcommand{\bfL}{\mathbf{L}}
\newcommand{\bfX}{\mathbf{X}}
\newcommand{\bfF}{\mathbf{F}}
\newcommand{\bfY}{\mathbf{Y}}
\newcommand{\bfZ}{\mathbf{Z}}
\newcommand{\bfM}{\mathbf{M}}
\newcommand{\bfP}{\mathbf{P}}
\newcommand{\bfQ}{\mathbf{Q}}
\newcommand{\bfT}{\mathbf{T}}
\newcommand{\bfV}{\mathbf{V}}
\newcommand{\bfW}{\mathbf{W}}
\newcommand{\bfI}{\mathbf{I}}
\newcommand{\bfR}{\mathbf{R}}
\newcommand{\bfS}{\mathbf{S}}
\newcommand{\bfU}{\mathbf{U}}
\newcommand{\bfO}{\mathbf{O}}
\newcommand{\bfzero}{\mathbf{0}}

\newcommand{\bfJ}{\mathbf{J}}

\newcommand{\rmT}{\mathrm{T}}

\newcommand{\calB}{\mathcal{B}}
\newcommand{\calD}{\mathcal{D}}
\newcommand{\calC}{\mathcal{C}}
\newcommand{\calT}{T}
\newcommand{\calF}{\mathcal{F}}
\newcommand{\calJ}{\mathcal{J}}
\newcommand{\calH}{\mathcal{H}}
\newcommand{\calL}{\mathcal{L}}
\newcommand{\calM}{\mathcal{M}}

\newcommand{\calK}{\mathcal{K}}
\newcommand{\calI}{\mathcal{I}}
\newcommand{\calQ}{\mathcal{Q}}
\newcommand{\calS}{\mathcal{S}}

\newcommand{\calW}{\mathcal{W}}
\newcommand{\calZ}{\mathcal{Z}}

\newcommand{\bfa}{\mathbf{a}}
\newcommand{\bfb}{\mathbf{b}}

\newcommand{\bfd}{\mathbf{d}}
\newcommand{\bfe}{\mathbf{e}}

\newcommand{\bfp}{\mathbf{p}}
\newcommand{\bfq}{\mathbf{q}}

\newcommand{\bfs}{\mathbf{s}}

\newcommand{\bfu}{\mathbf{u}}
\newcommand{\bfv}{\mathbf{v}}

\newcommand{\bfx}{\mathbf{x}}
\newcommand{\bfy}{\mathbf{y}}

\newcommand{\bt}{\begin{theorem}}
\newcommand{\et}{\end{theorem}}
\newcommand{\bl}{\begin{lemma}}
\newcommand{\el}{\end{lemma}}
\newcommand{\bp}{\begin{proposition}}
\newcommand{\ep}{\end{proposition}}
\newcommand{\bc}{\begin{corollary}}
\newcommand{\ec}{\end{corollary}}

\newcommand{\bd}{\begin{definition}\rm}
\newcommand{\ed}{\end{definition}}
\newcommand{\bex}{\begin{example}\rm}
\newcommand{\eex}{\end{example}}
\newcommand{\br}{\begin{remark}\rm}
\newcommand{\er}{\end{remark}}

\newcommand{\btbh}{\begin{table}[!ht]}
\newcommand{\etb}{\end{table}}
\newcommand{\bfgh}{\begin{figure}[!ht]}
\newcommand{\efg}{\end{figure}}

\newcommand{\bea}{\begin{eqnarray*}}
\newcommand{\eea}{\end{eqnarray*}}
\newcommand{\be}{\begin{eqnarray}}
\newcommand{\ee}{\end{eqnarray}}
\def\sspan{\mathop{\mathrm{span}}}
\def\rank{\mathop{\mathrm{rank}}}
\def\span{\mathop{\mathrm{span}}}

\newcommand{\Arg}{\mathop\mathrm{Arg}}

\makeatletter
\def\adots{\mathinner{\mkern2mu\raise\p@\hbox{.}
\mkern2mu\raise4\p@\hbox{.}\mkern1mu
\raise7\p@\vbox{\kern7\p@\hbox{.}}\mkern1mu}}
\newcommand{\l@abcd}[2]{\hbox to\textwidth{#1\dotfill #2}}
\makeatother

\def\unit{\mathfrak{i}}

\def\code#1{\texttt{#1}}

\begin{document}

\title{Low-rank signal subspace: \\parameterization, projection and signal estimation
\footnote{The reported study was funded by RFBR, project number 20-01-00067}}

\author{Nikita Zvonarev\footnote{{Faculty of Mathematics and Mechanics}, {St.Petersburg State University}, {Universitetskaya nab. 7/9, St.Petersburg, 199034}, {Russia}, nikitazvonarev@gmail.com}, Nina Golyandina\footnote{{Faculty of Mathematics and Mechanics}, {St.Petersburg State University}, {Universitetskaya nab. 7/9, St.Petersburg, 199034}, {Russia}, n.golyandina@spbu.ru, nina@gistatgroup.com}}

\maketitle

\abstract{\small The paper contains several theoretical results related to the weighted nonlinear least-squares problem for low-rank signal estimation, which can be considered as a Hankel structured low-rank approximation problem. A parameterization of the subspace of low-rank time series connected with generalized linear recurrence relations (GLRRs) is described and its features are investigated. It is shown how the obtained results help to describe the tangent plane, prove optimization problem features and construct stable algorithms for solving low-rank approximation problems. For the latter, a stable algorithm for constructing the projection onto a subspace of time series that satisfy a given GLRR is proposed and justified. This algorithm is used for a new  implementation of the known Gauss-Newton method using the variable projection approach.  The comparison by stability and computational cost is performed theoretically and with the help of an example.
}

\tableofcontents

\section{Introduction}\label{sec:intro}
Consider a class of time series, which corresponds to a model of signals in many real-life problems.
Denote $\tsS = (s_1, \ldots, s_N)^\rmT$ a signal of length $N$.
       The rank of $\tsS$  is defined as follows.
    For a given integer $L$, $1<L<N$, called the window length, we define the embedding operator $\calT_{L}:\; \spR^{N} \to \spR^{L \times (N-L+1)}$,
    which maps $\tsS$ into a Hankel $L\times (N-L+1)$ matrix, by
    	\begin{equation}
\label{eq:embedding}
	\calT_{L}(\tsS) = \begin{pmatrix}
	s_1 & s_2 & \hdots & s_{N-L+1} \\
	s_2 & s_3 & \hdots & \vdots \\
	\vdots & \vdots & \hdots & s_{N-1} \\
	s_{L} & s_{L+1} & \hdots & s_{N}
	\end{pmatrix}.
	\end{equation}
The columns of $\calT_{L}(\tsS)$ are sequential lagged vectors; this is why $\calT_{L}(\tsS)$ is often called the trajectory matrix of $\tsS$.
	We say that the signal $\tsS$ has rank $r < N/2$ if $\rank \calT_{r+1}(\tsS) = r$.
It is known that $\rank \calT_{r+1}(\tsS) = r$ if and only if $\rank \calT_{L}(\tsS) = r$ for any $L$ such that $\min(L, N-L+1) > r$ (see \cite[Corollary 5.1]{Heinig1984} for the proof).

 Denote $\calD_r$ the set of series of rank $r$.
	Since the set $\calD_r$ is not closed, we will also consider its closure $\overline{\calD_r}$. It is well-known that $\overline{\calD_r}$ consists of series of rank not larger than $r$ (this result can be found in \cite[Remark 1.46]{iarrobino1999power} for the complex case; the real-valued case is considered in Section~\ref{sec:closure}).

For a sufficiently large time series length $N$, it is well known that any signal in the form
   \begin{equation}
    \label{eq:model}
    s_n = \sum_{k=1}^d P_{m_k}(n)\exp(\alpha_k n) \sin(2\pi \omega_k n + \phi_k),
    \end{equation}
 where   $P_{m_k}(n)$ are polynomials in $n$ of degree $m_k$,
 has rank $r$, which is determined by the parameters $m_k$, $\alpha_k$ and $\omega_k$
(see Section~\ref{sec:rank_calc} for explaining the correspondence between the form of \eqref{eq:model} and the rank $r$).
  In signal processing applications, the signal in the model  \eqref{eq:model} is usually a sum of sine waves \cite{Cadzow1988}
 or a sum of damped sinusoids \cite{Markovsky2008}.

In this study we consider  the `signal plus noise' model of time series:
$$
x_n=s_n+\epsilon_n, \;\;\; n=1,2, \ldots, N.
$$
Denote by
$\tsX = (x_1, \ldots, x_N)^\rmT$ , $\tsS = (s_1, \ldots, s_N)^\rmT$  and $\bm\epsilon = (\epsilon_1, \ldots, \epsilon_N)^\rmT $
the vectors of observations,  signal values and errors respectively.

One of approaches for estimating the unknown signal values $s_n$ is solving the weight least-squares (WLS) problem
	\begin{equation}
    \label{eq:wls}
	\tsY^\star = \argmin_{\tsY \in \overline{\calD_r}} \| \tsX - \tsY \|_{\bfW},
	\end{equation}
   where $\bfW$ is a weight matrix and $\|\tsZ\|_{\bfW}^2 = \tsZ^\rmT \bfW \tsZ$.
 If noise $\bm\epsilon$ is Gaussian with covariance matrix $\bm\Sigma$ and zero mean, the WLS estimate with the weight matrix $\bfW = \bm\Sigma^{-1}$ is the maximum likelihood estimator (MLE). The same is true if the covariance matrix is scaled by a constant.

It is discussed in Section~\ref{sec:lrr} that each time series $\tsS$ from $\overline\calD_r$ is characterized by a vector $\bfa\in \spR^{r+1}$, which provides the coefficients of a generalized linear recurrence relation (GLRR) governing the time series, i.e. $\bfa^\rmT \calT_{r+1}(\tsS)$ is the zero vector. For each $\bfa$, we can consider the space $\calZ(\bfa)$ of signals governed by the GLRR with the given coefficients $\bfa$. Algorithms that use the variable projection method for solving the problem \eqref{eq:wls} include the projection onto $\calZ(\bfa)$ as a subproblem.

\textbf{The proposed approach.}
In this paper, the properties, which are based on the chosen family of local parameterizations of the low-rank time-series space related to generalized linear recurrent relations GLRR($\bfa$), are studied. We show that this parameterization is smooth and therefore allows one to consider different numerical optimization methods (e.g. Gauss-Newton method) for solving least-squares problems.
We prove (Theorem~\ref{th:tangent}) that the tangent subspace at the point $\tsS$, which is governed by a GLRR($\bfa$), can be described in terms of the GLRR($\bfa^2$).

The other contribution of the paper is the construction of a numerically stable algorithm of direct projection of a time series onto the space $\calZ(\bfa)$. This algorithm is fast for the case of a banded weight matrix $\bfW$, which corresponds to the case of autoregressive noise.

The proposed algorithm can be useful for numerical solutions of different approximation problems related to Hankel structured low-rank approximation (SLRA) problems; this is demonstrated by means of improving the stability of the known algorithm of low-rank time series approximation from \cite{Usevich2014}.

\textbf{Comments to general terminology}.
Let us explain the terminology, which we use. For each window length $L$, there is a one-to-one correspondence between a time series $\tsX_N$ of length $N$ and its $L$-trajectory matrix $\calT_L(\tsX_N)$ (we use this name taken from singular spectrum analysis) in $\spR^{L\times (N-L+1)}$. Different window lengths $L$ correspond to different (unweighted) matrix approximations. Therefore, the low-rank matrix approximations can be varied for different $L$. The notion of low-rank signals does not depend on $L$ and therefore is not related to matrices (generally speaking). Moreover, the solved problem \eqref{eq:wls} is stated in terms of time series, not in terms of matrices. For approximation by low-rank signals, $N$ weights are set for time series points, not for matrix entries. That is why we use the notion ``low-rank signals''.

Note that to have equivalent optimization problems for the matrix SLRA itself and for time series (vector) LRA, one should consider weighted versions and care about the correspondence of weights.
In \cite{Zvonarev.Golyandina2017} (see also a general description in \cite[Section 3.4]{Golyandina.etal2018}), the problem is solved as a matrix approximation problem with appropriate weights. In this paper, we consider the problem of time series low-rank approximation \eqref{eq:wls}.

\textbf{Structure of the paper}.
In  Section~\ref{sec:parameterization} we consider a parameterization of $\calD_r$ and its properties.
In Section~\ref{sec:optim} we apply the constructed parameterization for solving the optimization problem~\eqref{eq:wls} and introduce the known (VPGN) iterative method.
The algorithm VPGN is described in the way different from
that in \cite{Usevich2014}, since the description in \cite{Usevich2014} is performed for general SLRA problems
and therefore it is difficult to apply it to the particular case of Hankel SLRA for time series.
(For the convenience of readers, we include Table~\ref{table:defines} containing equivalent notations.)
In Sections~\ref{sec:ZofA} and~\ref{sec:proj_our} we propose the effective algorithm for implementing the projection onto $\calZ(\bfa)$.
In Section~\ref{sec:MGNand VPGN} we discuss the use of the proposed method of projection for stability improving the VPGN algorithm.
Section~\ref{sec:conclusion} concludes the paper.
Long proofs and technical details are relegated to the appendix.

\textbf{Main notation}.
In this paper, we use lowercase letters ($a$,$b$,\ldots) and also $L$, $K$, $M$, $N$ for scalars, bold lowercase letters ($\bfa$,$\bfb$,\ldots) for vectors, bold uppercase letters ($\bfA$,$\bfB$,\ldots) for matrices, and
the calligraphic font for sets. Formally, time series are vectors; however, we use the uppercase sans serif font ($\tsA$,$\tsB$,\ldots) for time series to distinguish them from ordinary vectors.
Additionally, $\bfI_{M}\in \spR^{M\times M}$ is the identity matrix, $\bm{0}_{M \times k}$ denotes the $M \times k$ zero matrix, $\bm{0}_M$ denotes the zero vector in $\spR^M$,
$\bfe_i$ is the $i$-th standard basis vector.

Denote $\bfb_{\vecrow{\calC}}$ the vector consisting of
the elements of a vector $\bfb$ with the numbers from a set $\calC$,
For matrices, denote $\bfB_{\row{\calC}}$ the matrix consisting
of rows of a matrix $\bfB$ with the numbers from $\calC$ and $\bfB_{\col{\calC}}$
the matrix consisting
of columns of a matrix $\bfB$ with the numbers from $\calC$.

Finally, we put a brief list of main common symbols and acronyms.\\
LRR is linear recurrence relation.\\
GLRR($\bfa$) is generalized LRR with the coefficients given by $\bfa$.\\
$\calD_r$ is the set of time series of rank $r$.\\
$\overline{\calD_r}$ is the set of time series of rank not larger than $r$.\\
$\calZ(\bfa) \in \spR^N$ is the set of time series of length $N$ governed by the minimal GLRR($\bfa$);
$\bfZ(\bfa)$ is the matrix consisting of its basis vectors.\\
$\calQ(\bfa)$ is the orthogonal complement to $\calZ(\bfa)$; $\bfQ(\bfa)$ is the matrix consisting of its special basis vectors in the form \eqref{op:Q}.\\
$\bfI_N$ is the $N\times N$ identity matrix.\\
$\bfW\in \spR^{N\times N}$ is a weight matrix.\\
$\winverse{\bfF}{\bfW}$ is the weighted pseudoinverse matrix; $\inverse{\bfF}$ stands for $\winverse{\bfF}{\bfI_N}$.\\
$\bfJ_{S}$ is the Jacobian matrix of a map $S$.\\
$\calT_M: \spR^N \rightarrow \spR^{M\times (N-M+1)}$ is the embedding operator, which constructs the $M$-trajectory matrix.\\
$\fullop$: $\spR^{r} \to \spR^{r+1}$ is the operator, which inserts $-1$ at the $\tau$ position.\\
$\Proj_{\calL,\bfW}$ is the $\bfW$-orthogonal projection onto $\calL$, $\Proj_{\bfL,\bfW}$ is the $\bfW$-orthogonal projection onto $\colspace(\bfL)$; if $\bfW$ is the identity matrix, it is omitted in the notation.\\
$S_{\tau}^\star(\Ai0) = \Proj_{\calZ(\fullop(\Ai0)), \bfW}(\tsX)$, where $\Ai0\in \spR^r$.

	\section{Parameterization of low-rank series}
    \label{sec:parameterization}
    \subsection{Generalized linear recurrence relations}
    \label{sec:lrr}

    It is well known \cite[Theorem 3.1.1]{Hall1998} that
    a time series of the form \eqref{eq:model} satisfies a linear recurrence relation (LRR) of some order $m$:
    \begin{equation}
    \label{eq:lrr}
    s_n = \sum_{k=1}^m b_k s_{n-k}, n = m+1, \ldots, N; b_m\neq 0.
    \end{equation}
    One time series can be governed by many different LRRs. The LRR of minimal order $r$ (it is unique) is called minimal.
    The corresponding time series has rank $r$. The minimal LRR uniquely defines the form of \eqref{eq:model} and the parameters $m_k$, $\alpha_k$, $\omega_k$.

    The relations \eqref{eq:lrr} can be expressed in vector form as
    $\bfa^\rmT \calT_{m+1}(\tsS) = \bfzero_{N-m}^\rmT$, where the vector $\bfa = (b_m, \ldots, b_1, -1)^\rmT \in \spR^{m+1}$.
    The vector $\bfa$ corresponding to the minimal LRR ($m=r+1$) and the first $r$  values of the series $\tsS$ uniquely determine the whole series  $\tsS$.
    Therefore, $r$ coefficients of an LRR of order $r$ and $r$ initial values
    ($2r$ parameters altogether) can be chosen as parameters of a time series of rank $r$.
    However, this parameterization does not describe the whole set $\calD_r$ \cite[Theorem 5.1]{Golyandina.etal2001}.

Let us generalize LRRs.
We say that a time series satisfies a generalized LRR (GLRR) of order $m$ if $\bfa^\rmT \calT_{m+1}(\tsS) = \bfzero_{N-m}^\rmT$ for some non-zero $\bfa \in \spR^{m+1}$;
we call this linear relation GLRR($\bfa$).
As well as for LRRs, the minimal GLRR can be introduced.
The difference between a GLRR  and an ordinary LRR is that the last coefficient in the GLRR is not necessarily non-zero and
therefore the GLRR does not necessarily set a recurrence.
However, at least one of the coefficients of the GLRR should be non-zero.
GLRRs correspond exactly to the first characteristic polynomial in \cite[Definition 5.4]{Heinig1984}.

Let us demonstrate the difference between LRR and GLRR by an example. Let $\tsS = (s_1,\ldots,s_N)^\rmT$ be a signal and $\bfa = (a_1, a_2, a_3)^\rmT$.
Then governing by GLRR($\bfa$) or LRR($\bfa$) means the same:
$a_1 s_i + a_2 s_{i+1} + a_3 s_{i+2} = 0$ for $i=1,\ldots,N-2$. For LRR($\bfa$), we state that $a_3  = -1$ (or just not equal to 0). Then this linear relation becomes a recurrence relation since $s_{i+2} = a_1 s_i + a_2 s_{i+1}$.
For GLRR($\bfa$), we assume that some of $a_i$ is not zero (or equal to $-1$). It may be $a_1$ or $a_2$ or $a_3$.

Any signal of rank $r$ satisfies a GLRR($\bfa$), where $\bfa\in \spR^{r+1}$. However, not each signal of rank $r$ corresponds to an LRR. E.g., $\tsS=(1,1,1,1,1,2)^\rmT$ has rank 2 and does not satisfy an LRR. However, it satisfies the GLRR($\bfa$) with $\bfa = (1,-1,0)^\rmT$. Therefore, we consider the parameterization with the help of GLRR($\bfa$). In fact, the same approach is used in \cite{ Usevich2014,Usevich2012}. It is indicated in Table~\ref{table:defines} that $\bfa$ in this paper corresponds to $R$ in \cite{ Usevich2014,Usevich2012}.

The following properties clarify the structure of the spaces $\calD_r$ and $\overline{\calD_r}$:
(a)
$\overline{\calD_r} = \{\tsY: \exists \bfa\in \spR^{r+1}, \bfa\neq \bfzero_{r+1}: \bfa^\rmT \calT_{r+1}(\tsS) = \bfzero_{N-r}^\rmT\}$ or, equivalently,
$\tsY \in \overline{\calD_r}$ if and only if there exists a GLRR($\bfa$) of order $r$, which governs $\tsY$;
(b)
$\tsY \in \calD_r$ if and only if there exists a GLRR($\bfa$) of order $r$, which governs $\tsY$, and this GLRR is minimal.

\subsection{Subspace approach}
\label{subsec:subspace_approach}
	Let  $\calZ(\bfa)$, $\bfa\in \spR^{r+1}$, be the space of time series of length $N$ governed by the GLRR($\bfa$);
that is, $\calZ(\bfa) = \{\tsS: \bfa^\rmT \calT_{r+1}(\tsS) = \bfzero_{N-r}^\rmT \}$.
Therefore $\overline{\calD_r} = \bigcup \limits_{\bfa} \calZ(\bfa)$.

Let  $\bfQ^{M, d}$ be the operator $\spR^{d+1} \to \spR^{M \times (M - d)}$, which is defined  by
\begin{equation}\label{op:Q}
\big(\bfQ^{M,d}(\bfb)\big)^\mathrm{T} = \begin{pmatrix}
b_1 & b_2 & \dots & \dots & b_{d+1} & 0 & \dots & 0 \\
0 & b_1 & b_2 & \dots &  \dots & b_{d+1} & \ddots & \vdots \\
\vdots & \ddots  & \ddots & \ddots & \ddots & \ddots & \ddots & 0 \\
0 & \dots & 0 & b_1 & b_2 & \ddots & \ddots & b_{d+1} \\
\end{pmatrix},
	\end{equation}
where $\bfb=(b_1, \ldots, b_{d+1})^\rmT \in \spR^{d+1}$.
Then the other convenient form of $\calZ(\bfa)$ is $$\calZ(\bfa) = \{\tsS: \bfQ^\rmT(\bfa) \tsS = \bfzero_{N-r} \},$$
where $\bfQ = \bfQ^{N,r}$. %; i.e., $\calZ(\bfa)$ is the left nullspace of $\bfQ(\bfa)$.

The following notation will be used below: $\calQ(\bfa) = \colspace(\bfQ(\bfa))$ and denote $\bfZ(\bfa)$ a matrix whose column vectors form a basis of $\calZ(\bfa)$. 

\subsection{Parameterization}
\label{sec:param}
Consider a series $\tsS_0\in \calD_r$, which satisfies a minimal GLRR($\bfa_0$) of order $r$ defined by a non-zero vector $\bfa_0 = (a_1^{(0)}, \ldots, $ $a_{r+1}^{(0)})^\rmT$.
Let us fix $\i0$ such that $a^{(0)}_{\i0} \neq 0$. Since GLRR($\bfa_0$) is invariant to multiplication by a constant,
we assume that $a^{(0)}_{\i0} = -1$. This condition on $\i0$ is considered to be valid hereinafter.
Let us build a parameterization of $\calD_r$ in the vicinity of $\tsS_0$; parameterization depends on the index $\i0$.
Note that we can not construct a global parameterization, since for different points of $\calD_r$ the index $\i0$,
which corresponds to a non-zero element of $\bfa_0$, can differ.

In the case of a series governed by an ordinary LRR($\bfa$), $\bfa\in \spR^{r+1}$, since the last coordinate of $\bfa$ is equal to $-1$,
the series is uniquely determined by the first $r$ elements of $\bfa$ and $r$  initial values of the series.
Then, applying the LRR to the initial data, which are taken from the series that is governed by the LRR,
we restore this series.

In the case of an arbitrary series from $\calD_r$, the approach is similar but a bit more complicated. For example,
we should take the boundary data ($\i0-1$ values at the beginning, and $r+1-\i0$ values at the end) instead of the $r$ initial values at the beginning of the series; also, the GLRR is not in fact recurrent (we keep the notation
to show that LRRs are a particular case of GLRRs).

Denote $\calI(\i0) = \{1,\ldots, N\} \setminus \{\i0,\ldots, N-r-1+\i0\}$ and
$\calK(\i0) = \{1,\ldots,r+1\} \setminus \{\i0\}$ two sets of size $r$.
The set $\calI(\i0)$ consists of the numbers of series values (we call them boundary data), which
are enough to find all the series values with the help of $\bfa$ (more precisely,
by elements of $\bfa$ with numbers from $\calK(\i0)$).
Then $\bfa_{\vecrow{\calK({\i0})} }\in \spR^{r}$
defines the vector consisting of the elements of a vector $\bfa \in \spR^{r+1}$
with the numbers from $\calK({\i0})$.

To simplify notation, let us introduce the operator $\fullop$: $\spR^{r} \to \spR^{r+1}$, which acts as follows. Let $\Ai0\in \spR^{r}$ and $\fullop(\Ai0) = \bfa$. Then $\bfa = (a_1, \ldots, $ $a_{r+1})^\rmT$ is such that $\bfa_{\vecrow{\calK({\i0})}} = \Ai0$ and $a_{\i0} = -1$; that is, $\Ai0 \in \spR^r$ is extended to $\bfa \in \spR^{r+1}$ by inserting $-1$ at the $\i0$-th position.
In this notation, $\bfa_{\vecrow{\calK({\i0})}} = \fullop^{-1}(\bfa)$.

Theorem~\ref{th:parameterization} defines the parameterization, which will be used in what follows.
The explicit form of this parameterization is given in Proposition~\ref{prop:parameterization}.
%From now, the same notation $\tsS$ is used for both the parameterizing mapping and for the series itself.

\begin{theorem} \label{th:parameterization}
	Let $\bfa_0\in \spR^{r+1}$, $a_\i0^{(0)} = -1$, and $\tsS_0 \in \calD_r$ satisfy the GLRR($\bfa_0$). Then there exists a unique one-to-one mapping $S_{\tau}: \spR^{2r} \to \calD_r$ between a neighborhood of the point $\left((\tsS_0)_{\calI(\i0)}, (\bfa_0)_{\calK(\i0)}\right)^\rmT \in \spR^{2r}$ and the intersection of a neighborhood of $\tsS_0$ with the set $\calD_r$, which satisfies the following relations:
for $\tsS = S_{\tau}(\Si0, \Ai0)$, where $\Si0, \Ai0 \in \spR^r$, we have
\begin{itemize}
\item
 $(\tsS)_{\vecrow{\calI(\i0)}} = \Si0$;
 \item
 $\tsS \in \calD_r$ is governed by the GLRR($\fullop(\Ai0)$).
 \end{itemize}
\end{theorem}

\begin{proposition} \label{prop:parameterization}
Let $\bfa_0\in \spR^{r+1}$, $a_\i0^{(0)} = -1$, and $\bfZ_0 \in \spR^{N \times r}$ consist of basis vectors of $\calZ(\bfa_0)$. Consider the parameterizing mapping $S_\tau$, introduced in Theorem~\ref{th:parameterization}.\\
	%\begin{enumerate}
	%	\item
1. Let $(\Si0, \Ai0)^\rmT \in \spR^{2r}$ and denote $\bfa = \fullop(\Ai0)$.
 Denote $\Proj_{\calZ(\bfa)}$ the orthogonal projection onto $\calZ(\bfa)$. Then for $\bfZ = \Proj_{\calZ(\bfa)} \bfZ_0\in \spR^{N\times r}$ and $\bfG = \bfZ \left(\bfZ_{\row{\calI({\i0})}}\right)^{-1}$, where $\bfZ_{\row{\calI({\i0})}}\in \spR^{r\times r}$, the mapping $S_\tau$ has the explicit form
		\begin{equation}\label{eq:param}
		\tsS = S_\tau(\Si0, \Ai0) = \bfG \Si0.
		\end{equation}
%		\item
2. The inverse of the mapping $S_\tau$ is given as follows. Let $\tsS = S_\tau(\Si0, \Ai0)$. Then
\begin{equation}
         \label{eq:param_rev}
		%\begin{aligned}
        \Si0 = (\tsS)_{\vecrow{\calI(\i0)}},\qquad
		\Ai0 = (-\hat{\bfa}/\hat{a}_\i0)_{\vecrow{\calK(\i0)}},
		%\end{aligned}
\end{equation}
		where $\hat{\bfa} = \hat{\bfa}(\tsS) = (\hat a_1, \ldots, \hat a_{r+1})^\rmT = \left(\bfI_{r+1} -  \Proj_{\calL(\tsS)}\right) \bfa_0$, $\calL(\tsS) = \colspace\left(\calT_{r+1}(\tsS)\right)$, $\Proj_{\calL(\tsS)}$ is the orthogonal projection onto $\calL(\tsS)$.
%	\end{enumerate}
\end{proposition}

\begin{proof}
See the proof of Theorem~\ref{th:parameterization} together with Proposition~\ref{prop:parameterization} in Section~\ref{sec:th:parameterization}.
\end{proof}

Note that for different series $\tsS_0\in \calD_r$ we have different parameterizations of $\calD_r$ in vicinities of $\tsS_0$.
Moreover, for a fixed $\tsS_0$, there is a variety of parameterizations provided by different choices of the index $\i0$.

\subsection{Smoothness of parameterization and derivatives}
\begin{theorem}
	\label{th:param_smooth}
Let $\bfa_0\in \spR^{r+1}$, $a_\i0^{(0)} = -1$, and $\tsS_0 \in \calD_r$ satisfy the GLRR($\bfa_0$). Then the parameterization $S_{\tau}(\Si0, \Ai0)$, which is introduced in Theorem~\ref{th:parameterization} and Proposition~\ref{prop:parameterization},
 is a smooth diffeomorphism
 between a neighborhood of the point $\left((\tsS_0)_{\calI(\i0)}, (\bfa_0)_{\calK(\i0)}\right)^\rmT \in \spR^{2r}$ and the intersection of a neighborhood of $\tsS_0$ with the set $\calD_r$.
\end{theorem}
    \begin{proof}
    	We need to show that $\Proj_{\calL(\tsS)}$ and $\Proj_{\calZ(\bfa)}$ from Proposition \ref{prop:parameterization} are smooth projections in the vicinity of $\bfS_0$ and $\bfZ_0$ respectively.
    	
    	Since $(\bfS_0)_{\col{\calJ}}$ has full rank, $\Proj_{\calL(\tsS)} = \bfS_{\col{\calJ}} \left(\left(\bfS_{\col{\calJ}}\right)^\rmT \bfS_{\col{\calJ}}\right)^{-1} \bfS^\rmT_{\col{\calJ}}$ is a smooth function in the vicinity of $\bfS_0$.
    	Since $\bfQ(\bfa)$ has full rank, see definition \eqref{op:Q}, $$\Proj_{\calZ(\bfa)} = \bfI_N - \bfQ(\bfa)\left(\bfQ^\rmT(\bfa) \bfQ(\bfa)\right)^{-1} \bfQ^\rmT(\bfa)$$ is smooth everywhere except $\bfa = \bfzero_{r+1}$.
    	
    	It is seen that the other mappings involved in the parameterization are smooth in the corresponding vicinities.
    \end{proof}

Let us consider the derivatives of the parameterizing mapping.
Let the series $\tsS$ belong to a sufficient small neighborhood of $\tsS_0$ and be parameterized as $\tsS = S_{\tau}(\Si0, \Ai0)$.
Denote $\bfJ_{S_{\tau}} = \bfJ_{S_{\tau}}(\Si0, \Ai0) \in \spR^{N\times 2r}$ the Jacobian matrix of  $S_{\tau}(\Si0, \Ai0)$.

By definition, the tangent subspace at the point $\tsS$ coincides with $\colspace\left(\bfJ_{S_{\tau}}(\Si0, \Ai0)\right)$.
Note that the tangent subspace is invariant with respect to the choice of a certain parameterization of $\calD_r$ in the vicinity of $\tsS$.

Define by $\bfa^2$ the acyclic convolution of $\bfa$ with itself:
\begin{equation*}
	\bfa^2 = (a^{(2)}_i) \in \spR^{2r+1}, \quad a^{(2)}_i = \sum_{j=\max(1, i - r)}^{\min(i, r+1)} a_j a_{i - j + 1}.
\end{equation*}
		
\begin{theorem}
\label{th:tangent}
	The tangent subspace to $\calD_r$ at the point $\tsS$ has dimension $2r$ and is equal to $\calZ(\bfa^2)$.
\end{theorem}

Let us start the proof with two lemmas.
    It is convenient to separate the parameters ($2r$ arguments of the mapping $S_\tau$) into two parts, $\Si0$ and $\Ai0$.
    Then $\bfJ_{S_\tau} = [ \bfF_\bfs:\bfF_\bfa ]$, where $\bfF_\bfs = (\bfJ_{S_\tau})_{\col{\{1, \dots, r\}}}$, $\bfF_\bfa = (\bfJ_{S_\tau})_{\col{\{r+1, \dots, 2r\}}}$. Let $\bfa = \fullop(\Ai0)$.

    \begin{lemma}
    	\label{eqa:derivS}
    	$\bfQ^\rmT(\bfa) \bfF_\bfs = \bm{0}_{(N-r) \times r}$; $\colspace(\bfF_\bfs) = \calZ(\bfa)$.
    \end{lemma}	
    \begin{proof}
    	Let $\bfF_\bfs = [F_{s,1}:\ldots:F_{s,r}]$. Consider the equality $\bfQ^\rmT(\bfa) S_\tau(\Si0, \Ai0) = \bfzero_{N-r}$ and differentiate it with respect to $(\Si0)_{\vecrow{(i)}}$. We obtain $\bfQ^\rmT(\bfa) F_{s, i} = \bfzero_{N-r}$, which means that $\colspace(\bfF_\bfs) \subset \calZ(\bfa)$.
    	The fact $(\bfF_{\bfs})_{\row{\calI({\i0})}} = \bfI_r$ completes the proof.
    \end{proof}

    \begin{lemma}
    	\label{eqa:derivA}
    	$\bfQ^\rmT(\bfa) \bfF_{\bfa} = \bfM$, where $\bfM = - (\bfS_{\row{\calK({\i0})}})^\rmT$ and $\bfS = \calT_{r+1}(\tsS)$;
    	$\colspace(\bfF_\bfa) \subset \calZ(\bfa^2)$.
    \end{lemma}	
    \begin{proof}
    	Let $\bfF_\bfa = [F_{a,1}:\ldots:F_{a,r}]$. Consider the equality $\bfQ^\rmT(\bfa) S_\tau(\Si0, \Ai0) = \bfzero_{N-r}$ and differentiate it with respect to $(\Ai0)_{\vecrow{(i)}}$, i.e. $i$-th element of $\bfa_{\calK({\i0})} = \Ai0$, $i = 1, \ldots, r$. Then we obtain $\bfQ^\rmT(\bfe_{j}) \tsS + \bfQ^\rmT(\bfa) F_{a, i} = \bfzero_{N-r}$, where $\bfe_j\in \spR^{r+1}$ and $j=(\calK({\i0}))_i$ is $i$-th element of $\calK({\i0})$. (Note that $\bfQ^\rmT(\bfe_{j}) \tsS$ is the $j$-th column of the transposed $(r+1)$-trajectory matrix $\bfS^\rmT$.) Therefore, the equation $\bfQ^\rmT(\bfa) \bfF_{\bfa} = - (\bfS_{\row{\calK({\i0})}})^\rmT$ is proved.
    	
    	To prove the second statement of the lemma, let us take the matrix $\bfQ^{N-r, r}(\bfa) \in \spR^{N \times (N - r)}$. Due to the first statement, the equality $\left(\bfQ^{N-r, r}(\bfa)\right)^\rmT \left(\bfQ^{N, r}(\bfa)\right)^\rmT \bfF_{\bfa} = \bm{0}_{(N-2r) \times r}$ is valid. From \cite[Sections 2.1 and 2.2]{Usevich2017} it follows that then we have $\left(\bfQ^{N-r, r}(\bfa)\right)^\rmT \left(\bfQ^{N, r}(\bfa)\right)^\rmT = \bfQ^\rmT(\bfa^2)$. Therefore, $\bfQ^\rmT(\bfa^2) \bfF_{\bfa} = \bm{0}_{(N-2r) \times r}$.
    \end{proof}

    Now we can prove Theorem \ref{th:tangent}.
    \begin{proof}
    	It follows from Lemma \ref{eqa:derivS} that
    	\begin{equation*}
    	\bfQ^\rmT(\bfa^2) \bfF_{S} = (\bfQ^{N-r, r}(\bfa))^\rmT \bfQ^\rmT(\bfa) \bfF_{S} = \bm{0}_{(N-2r) \times r}.
    	\end{equation*}
    	Therefore,  $\colspace(\bfJ_{S_\tau}) \subset \calZ(\bfa^2)$. Also, $\tsS \in \calZ(\bfa^2)$. Since we have a diffeomorphism at the point $\tsS$, the Jacobian matrix $\bfJ_{S_\tau}$ has full rank $2r$. Hence, $\colspace(\bfJ_{S_\tau}) = \calZ(\bfa^2)$.
    \end{proof}

\section{Parameterization and the low-rank optimization problem}
\label{sec:optim}
Let us consider the problem \eqref{eq:wls}.
First, note that we  search for a local minimum. Then, since the objective function is smooth in the considered
parameterization, one can apply the conventional weighted version of the Gauss-Newton method (GN), see \cite{nocedal2006numerical} for details.
However, this approach appears to be numerically unstable and has a high computational cost.

In \cite{Usevich2014}, the variable-projection method (VP) is used
for solving the minimization problem. When the reduced minimization problem
is solved again by the Gauss-Newton method; we will refer to it as VPGN.

Note that the considered methods are used for solving a weighted least-squares problem and therefore
we consider their weighted versions, omitting `weighted' in the names of the methods.

Let us introduce notation, which is used in this section.
For some matrix $\bfF = \spR^{N\times p}$, define its weighted pseudoinverse \cite{Stewart1989}
$\winverse{\bfF}{\bfW} = (\bfF^\rmT \bfW \bfF)^{-1}\bfF^\rmT \bfW$;
this pseudoinverse arises in the solution of the linear weighted least-squares problem $\min_\bfp \|\bfy - \bfF \bfp\|_{\bfW}^2$
with $\bfy \in \spR^{N}$, since its solution is equal to $\bfp_{\mathrm{min}} = \winverse{\bfF}{\bfW} \bfy$. In the particular case $\bfW = \bfI_N$, $\winverse{\bfF}{\bfW}$ is the ordinary pseudoinverse; we will denote it $\inverse{\bfF}$.
Denote the projection (it is oblique if $\bfW$ is not the identity matrix) onto the column space $\calF$ of a matrix $\bfF$ as $\Proj_{\bfF, \bfW} = \bfF \winverse{\bfF}{\bfW}$.
If it is not important which particular basis of $\calF$ is considered, we use the notation $\Proj_{\calF, \bfW}$.

\begin{remark}
\label{rem:complexF}
If the matrix $\bfF$ is complex, the above formulas and considerations are still valid with the change of the transpose $\bfF^\rmT$ to the complex conjugate $\bfF^*$.
\end{remark}

\subsection{Properties of the optimization problem}
 The following lemma shows that the global minimum of \eqref{eq:wls} belongs to $\calD_r$ for the majority of $\tsX$. Therefore, it is sufficient
 to find the minimum in the set of series of exact rank $r$.

 \begin{lemma}\label{lemma:minindr}
	Let $\tsX \notin \overline{\calD_r} \setminus \calD_r$ and $\bfW$ be positive definite. Then any point of the global minima in the problem \eqref{eq:wls} belongs to $\calD_r$.
\end{lemma}
    \begin{proof}
    	Assume the contrary. Denote $\tsS^\star = \tsS_0$ a point of global minimum in the problem \eqref{eq:wls} and assume that $\tsS_0 \in \calD_{r_0}$, $r_0 < r$, is such that $\tsS_0$ satisfies a GLRR($\bfa_0$), $\bfa_0 = (a_1, \ldots, a_{r_0+1})^\rmT \in \spR^{r_0 + 1}$. Construct $N$ linearly independent exponential series $\tsS^{(i)}$ of length $N$, $\tsS^{(i)} = (e^{\lambda_i}, e^{2 \lambda_i}, \ldots, e^{N \lambda_i})^\rmT$, which are governed by the GLRR($\bfa^{(i)}$) with $\bfa^{(i)} = (e^{\lambda_i}, -1)$, $i=1,\ldots,N$, where all $\lambda_i$ are different. Then for any real $\alpha$ we have
    	$\tsS_0 + \alpha \tsS^{(i)} \in \overline{\calD_r}$ since the series $\tsS_0 + \alpha \tsS^{(i)}$ is governed by the GLRR($\bfb_i$) with $\bfb_i = (e^{\lambda_i} a_1, e^{\lambda_i} a_2 - a_1, e^{\lambda_i} a_3 - a_2, \ldots, e^{\lambda_i} a_{r_0+1} - a_{r_0}, -a_{r_0 + 1})^\rmT \in \spR^{r_0 + 2}$.
    	
    	Denote $\langle \tsZ, \tsY \rangle_{\Sigminus} = \tsZ^\rmT \Sigminus \tsY$ the weighted inner product  corresponding to the norm $\| \cdot \|_{\Sigminus}$. By the condition of the lemma, $\tsX - \tsS_0 \ne \bm{0}_N$. Consider the inner products $\langle \tsX - \tsS_0, \tsS^{(i)} \rangle_{\Sigminus}$, $i = 1, 2, \ldots, N$. Since $\tsS^{(i)}$, $i=1,\ldots,N$, form a basis of $\spR^N$, there exists an index $j$ such that $\langle \tsX - \tsS_0, \tsS^{(j)} \rangle_{\Sigminus} \ne 0$. Let us take $\tsS_1 = \tsS_0 + \frac{\langle \tsX - \tsS_0, \tsS^{(j)} \rangle_{\Sigminus}}{\langle \tsS^{(j)}, \tsS^{(j)} \rangle_{\Sigminus}} \tsS^{(j)}$ governed by the GLRR($\bfb_i$) (hence, $\tsS_1$ belongs to $\overline{\calD_r}$), and show that $\|\tsX - \tsS_1\|_{\Sigminus} < \|\tsX - \tsS_0\|_{\Sigminus}$. Indeed,
    	\begin{equation*}
    	\langle \tsX - \tsS_0, \tsX - \tsS_0  \rangle_{\Sigminus} - \langle \tsX - \tsS_1, \tsX - \tsS_1  \rangle_{\Sigminus} = \frac{\left( \langle \tsX - \tsS_0, \tsS^{(j)} \rangle_{\Sigminus} \right)^2}{\langle \tsS^{(j)}, \tsS^{(j)} \rangle_{\Sigminus}}>0.
    	\end{equation*}
    	We obtain a contradiction to the initial assumption that $\tsS_0 = \tsS^\star$ is a point of global minimum in the problem \eqref{eq:wls}.
    \end{proof}

Thus, the problem \eqref{eq:wls} can be considered as a minimization problem in $\calD_r$; therefore, in the chosen parameterization of $\calD_r$ (see Section~\ref{sec:param}),
the problem \eqref{eq:wls} in the vicinity of $\tsS_0$ has the form
	\begin{equation}
    \label{eq:wlsP}
	\bfp^\star = \argmin_{\bfp} \| \tsX - S(\bfp) \|_{\bfW},
	\end{equation}
where $\bfp=(\Si0, \Ai0)$, $S = S_{\tau}$.
Since $S(\bfp)$ is a differentiable function of $\bfp$ due to Theorem~\ref{th:param_smooth}
for an appropriate choice of $\i0$, numerical methods like the Gauss-Newton method can be applied to the solution of \eqref{eq:wlsP}.

The following theorem helps to detect if the found solution is a local minimum. Recall that $\calZ(\bfa^2)$ determines the
tangent subspace (Theorem~\ref{th:tangent}).

\begin{lemma}[Necessary conditions for local minima]\label{lemma:locminnec}
	Let $\bfW$ be positive definite. If the series $\tsX_0 \in \calD_r$ which is governed by a GLRR($\bfa_0$) provides a local minimum in the problem \eqref{eq:wls},
then $\Proj_{\calZ(\bfa_0^2), \Sigminus}(\tsX - \tsX_0) = \bfzero_N$.
\end{lemma}
\begin{proof}
	Let us take an appropriate index $\i0$ together with the parameterization $S_{\tau}(\Si0, \Ai0)$ introduced in Theorem \ref{th:parameterization}.
Due to Theorem \ref{th:param_smooth}, the objective function $\|\tsX -  S_{\tau}(\Si0, \Ai0)\|^2_\bfW$ is smooth in the vicinity of $\left((\bfs_0)_{\calI(\i0)}, (\bfa_0)_{\calK(\i0)}\right)^\rmT \in \spR^{2r}$. Theorem~\ref{th:tangent} together with \cite[Theorem 2.2]{nocedal2006numerical}, which formulates the necessary conditions for a minimum in a general case.
\end{proof}

Note that Lemma~\ref{lemma:locminnec} provides the necessary condition only. According to \cite[Theorem 2.3]{nocedal2006numerical},
 sufficient conditions for a minimum include positive definiteness of the Hessian of the objective function.
 For the present, we cannot theoretically check this positive definiteness.

\subsection{The Gauss-Newton method with variable projection (VPGN)}
\label{sec:VPGN}

A variation from the standard way of the use of iterative methods is that the parameterization $S_{\tau}(\bfp)$, $\bfp=(\Si0, \Ai0)$ (which is based on $\i0$)
is changed at each iteration in a particular way.
At $(k+1)$-th iteration, the parameterization is constructed in the vicinity of  $\bfa_0 = \bfa^{(k)}$. The index $\i0$, which determines the parameterization, is chosen in such a way to satisfy $a^{(0)}_\i0 \neq 0$.
We propose the following approach to the choice of $\i0$.
Let $\i0$ be the index of the maximum absolute entry of $\bfa_0$. Since the parameterization is invariant to the multiplication of $\bfa_0$ by a constant, it can be assumed that $a^{(0)}_\i0 = -1$ and $|a^{(0)}_i| \le 1$ for any $i$, $1 \le i \le r+1$.

The explicit form of the parameterization $S_{\tau}(\Si0, \Ai0) = \bfG(\Ai0) \Si0$ given in \eqref{eq:param}, where $\Si0$ is presented in $S_{\tau}(\Si0, \Ai0)$ in a linear manner,
allows one to apply the variable projection principle (see \cite{Golub.Pereyr2003} for the case of the Euclidean norm).
Thus, the parameter $\Si0$ can be eliminated and the problem \eqref{eq:wls} is reduced to
\begin{equation}\label{eq:wlsvp_set}
\tsY^\star = \argmin_{\tsY \in \calD_r^\star} \| \tsX - \tsY \|_{\bfW},
\end{equation}
where $\calD_r^\star = \{\Proj_{\calZ(\fullop(\Ai0)), \bfW}(\tsX) \mid \Ai0 \in \spR^{r}\} \subset \overline \calD_r$.

Denote
\begin{equation}
S_{\tau}^\star(\Ai0) = \Proj_{\calZ(\fullop(\Ai0)), \bfW}(\tsX).
\end{equation}

Therefore, we can present the problem \eqref{eq:wlsvp_set} in terms of the parameter $\Ai0$ only:
\begin{equation}\label{eq:wlsvp1}
\Ai0^\star = \argmin_{\Ai0 \in \spR^r} \| \tsX - S_{\tau}^\star(\Ai0) \|_{\bfW},
\end{equation}
Thus, for the numerical solution of the equation \eqref{eq:wls}, it is sufficient to consider iterations for the nonlinear part of the parameters.
This is the VP approach used in \cite{Usevich2014, Usevich2012}.

Let us denote $\bfJ_{S_{\tau}^\star}(\Ai0)$ the Jacobian matrix of $S_{\tau}^\star(\Ai0)$.
Then the iterations of the Gauss-Newton method for solving the problem \eqref{eq:wlsvp1} have the form
\begin{equation}
\label{eq:gauss_simple}
\Ai0^{(k+1)} = \Ai0^{(k)} + \gamma \winverse{\bfJ_{S_{\tau}^\star}(\Ai0^{(k)})}{\bfW} (\tsX - S_{\tau}^\star(\Ai0^{(k)})).
\end{equation}

An explicit form of $\bfJ_{S_{\tau}^\star}(\Ai0^{(k)})$ is described below.

\subsubsection{Formulas for calculating the iteration step in VPGN} \label{sec:MUdetails}
An explicit form of the step \eqref{eq:gauss_simple} is contained in \cite[Proposition 3]{Usevich2014}.
Here we write down the formulas in our notation and also present a new form for the Jacobian $\bfJ_{\tsS_\tau^\star}$, which is more convenient for implementation.
\begin{lemma}
	\label{th:varproj}
	Let $\bfW$ be positive definite. The projection $\Proj_{\calZ(\bfa), \bfW}$ can be calculated as
	\begin{equation} \label{eq:spZaa_kostya}
	\Proj_{\calZ(\bfa), \bfW} \tsX = \left( \bfI_N - \bfW^{-1} \bfQ(\bfa) \bm\Gamma^{-1}(\bfa) \bfQ^\rmT(\bfa) \right) \tsX,
	\end{equation}
	where $\bm\Gamma(\bfa) = \bfQ^\rmT(\bfa) \Sigminus^{-1} \bfQ(\bfa)$.\\
	The columns of $\bfJ_{\tsS_\tau^\star}$ has the form
	\begin{equation} 	\label{eq:vpformula}
	(\bfJ_{\tsS_\tau^\star})_{\col{i}} = -\bfW^{-1} \bfQ(\bfa) \bm\Gamma^{-1}(\bfa) \bfQ^\rmT(\bfe_j) \Proj_{\calZ(\bfa), \bfW} \tsX
    -\Proj_{\calZ(\bfa), \bfW} \bfW^{-1} \bfQ(\bfe_j) \bm\Gamma^{-1}(\bfa) \bfQ^\rmT(\bfa) \tsX,
	\end{equation}
	where $\bfa = \fullop(\Ai0)$ and $j = (\calK({\i0}))_i$ is $i$-th element of $\calK({\i0})$.
	
\end{lemma}
\begin{proof}
    The equality
	\begin{equation*}
	S_\tau^\star(\Ai0) = \Proj_{\calZ(\fullop(\Ai0)), \bfW}(\tsX)
	\end{equation*}
	corresponds to the solution of the following quadratic problem:
	\begin{equation}
    \label{eq:lin_constraint}
	S_\tau^\star(\Ai0) = \argmin_{\substack{\tsY: \  \bfQ^\rmT(\bfa) \tsY=0}} \left(\frac{1}{2} \tsY^\rmT \bfW\tsY - \tsY \bfW \tsX \right).
	\end{equation}
    The problem \eqref{eq:lin_constraint} is the equality-constrained quadratic optimization problem, which can be written as a linear system \cite[Section 16.1]{nocedal2006numerical}. The Schur-complement method described in \cite[Section 16.2]{nocedal2006numerical} provides the expression \eqref{eq:spZaa_kostya} after substituting the corresponding notation.
	
	Proof of equality \eqref{eq:vpformula} is done by taking derivatives of \eqref{eq:spZaa_kostya} with respect to $a_j$:
	\begin{multline*}
	(\Proj_{\calZ(\bfa), \bfW} \tsX)'_{a_j} = - \bfW^{-1} \bfQ(\bfe_j) \left(\bfQ^\rmT(\bfa) \bfW^{-1} \bfQ(\bfa)\right)^{-1} \bfQ^\rmT(\bfa) \tsX \\
	-\bfW^{-1} \bfQ(\bfa) \left(\bfQ^\rmT(\bfa) \bfW^{-1} \bfQ(\bfa)\right)^{-1} \bfQ^\rmT(\bfe_j) \tsX \\+  \bfW^{-1} \bfQ(\bfa) \left(\bfQ^\rmT(\bfa) \bfW^{-1} \bfQ(\bfa)\right)^{-1}  \\
\times \left(\bfQ^\rmT(\bfe_j) \bfW^{-1} \bfQ(\bfa) + \bfQ^\rmT(\bfa) \bfW^{-1} \bfQ(\bfe_j) \right) \left(\bfQ^\rmT(\bfa) \bfW^{-1} \bfQ(\bfa)\right)^{-1} \\ \times \bfQ^\rmT(\bfa) \tsX
	=-\bfW^{-1} \bfQ(\bfa) \left(\bfQ^\rmT(\bfa) \bfW^{-1} \bfQ(\bfa)\right)^{-1} \bfQ^\rmT(\bfe_j) \Proj_{\calZ(\bfa), \bfW} \tsX \\
    -\Proj_{\calZ(\bfa), \bfW} \bfW^{-1} \bfQ(\bfe_j) \left(\bfQ^\rmT(\bfa) \bfW^{-1} \bfQ(\bfa)\right)^{-1} \bfQ^\rmT(\bfa) \tsX.
	\end{multline*}
\end{proof}

\section{Calculation of a particular orthonormal basis of $\calZ(\bfa)$}
\label{sec:ZofA}
\label{fouriermethod}
 In this section, we consider the construction of
such orthonormal bases that allow one to calculate the projections to $\calZ(\bfa)$ with improved precision.
The constructed algorithms can also be used to improve the numerical stability of the iteration step \eqref{eq:gauss_simple} of the VPGN method.

\subsection{Circulant matrices and construction of the basis}
{Denote $\bfZ(\bfa)$ a matrix consisting of basis vectors of the subspace $\calZ(\bfa)$ for some coefficient vector $\bfa$.}
Despite the series are real-valued, we construct a complex-valued basis of the complexification of $\calZ(\bfa)$, since this does not affect the result of the projection $\Proj_{\calZ(\bfa), \bfW} \bfv$ for any real
vector $\bfv$ and real matrix $\bfW$. Thus, we want to find a matrix $\bfZ(\bfa) = \bfZ \in \spC^{N \times r}$ of full rank to satisfy
%\be
%\label{eq:lineqMGN}
$
\bfQ^\rmT(\bfa) \bfZ = \bm{0}_{(N-r) \times r}.
$
%\ee

The matrix $\bfQ^\rmT(\bfa)$ is a partial circulant. Let us extend $\bfQ^\rmT(\bfa)$ to the circulant matrix $\bfC(\bfa)$ of $\bfa \in \spR^{r+1}$:
\begin{equation}\label{eq:circulant}
\bfC(\bfa) = \begin{pmatrix}%\left(\begin{smallmatrix}
a_1 & a_2 & \dots & \dots & a_{r+1} & 0 & \dots & 0 \\
0 & a_1 & a_2 & \dots &  \dots & a_{r+1} & \ddots & \vdots \\
\vdots & \ddots  & \ddots & \ddots & \ddots & \ddots & \ddots & 0 \\
0 & \dots & 0 & a_1 & a_2 & \ddots & \ddots & a_{r+1} \\
a_{r+1} & \ddots & \ddots & \ddots & \ddots & \ddots & \ddots & \vdots \\
\vdots & \ddots & \ddots & \ddots & \ddots & \ddots & \ddots & \vdots \\
a_3& \dots & a_{r+1} & 0 & \dots & 0 & a_1 & a_2 \\
a_2& \dots & \dots & a_{r+1} & 0 & \dots  & 0 & a_1
\end{pmatrix}.%\end{smallmatrix}\right).
\end{equation}
Then $\bfv \in \calZ(\bfa)$ if and only if  $\bfC (\bfa) \bfv \in \span(\bfe_{N-r+1}, \ldots, \bfe_N)$, $\bfe_i\in \spR^{r+1}$.
If $\bfC(\bfa)$ has full rank, then we can find the basis vectors $\bfv_k$ solving the systems of linear equations
\be
\label{eq:lineqMGN}
    \bfC (\bfa) \bfv_k = \bfe_{N-k+1},\ k=1,\ldots,r,
\ee
with the computational cost of the order $O(r N \log N)$, since the calculations can be performed with the help of the discrete Fourier transform \cite{Davis2012} by fast Fourier transform (FFT), and then applying orthonormalization to the columns of $\bfV_r = \left[ \bfv_1 : \ldots :\bfv_r \right]$.

Denote $\calF_N$ and $\calF^{-1}_N$ the Fourier transform and the inverse Fourier transform for series of length $N$, respectively. That is,
for $\bfx = (x_0, \ldots, $ $x_{N-1})^\rmT \in \spC^N$ we have $\calF_N(\bfx) = \bfy = (y_0, \ldots, y_{N-1})^\rmT \in \spC^N$, where
$y_k = \frac{1}{\sqrt{N}} \sum_{j = 0}^{N-1} x_j \exp\big(-\frac{\unit 2 \pi k j}{N}\big)$.
Define $\calF_N(\bfX) = [\calF_N(\bfx_1): \ldots: \calF_N(\bfx_r)]$, where $\bfX=[\bfx_1 : \ldots : \bfx_r]$; the same for $\calF_N^{-1}(\bfY)$.

 Let
 \be
 \label{eq:pol_z}
 g_{\bfa}(z) = \sum_{k=0}^{r} a_{k+1}z^k
 \ee
  be the complex polynomial with coefficients $\bfa = (a_1,\ldots,a_{r+1})^\mathrm{T}$; we do not assume that the leading coefficient is non-zero.

The following lemma is a direct application of the theorem about the solution of a linear system of equations given by a circulant matrix \cite{Davis2012}.

\begin{lemma}
	\label{lemma:ev_circulant_pre}
Denote $\bfV_r = \calF_N^{-1}(\bfA_g^{-1} \bfR_r)$, where $\bfA_g = \diag((g_\bfa(\omega_0), \ldots, g_\bfa(\omega_{N-1}))^\rmT)$ for $\omega_j = \exp\big(\frac{\unit 2 \pi  j}{N}\big)$ and $\bfR_r = \calF_N([ \bfe_{N-r}:\ldots:\bfe_{N}])$.
Then $\bfQ^\rmT(\bfa) \bfV_r = \bm{0}_{(N-r) \times r}$, that is, $\colspace(\bfV_r) = \calZ(\bfa)$.
Herewith, the diagonal of the matrix $\bfA_g$ consists of the eigenvalues of the circulant matrix $\bfC(\bfa)$.\\ \end{lemma}

\begin{remark} \label{rem:fourier_orthogon}
	%\begin{enumerate}
		%\item
1. Let $\bfZ = \mathrm{orthonorm}(\bfV_r)$ be a matrix consisting of orthonormalized columns of the matrix $\bfV_r = \calF_N^{-1}(\bfA_g^{-1} \bfR_r)$ given in Lemma \ref{lemma:ev_circulant_pre}. Then $\bfZ$ is a matrix whose columns form an orthonormal basis of $\calZ(\bfa)$. Indeed, since $\bfQ^\rmT(\bfa)\bfV_r = \bfzero_{(N-r)\times r}$, we have $\bfQ^\rmT(\bfa)\bfZ = \bfzero_{(N-r)\times r}$.\\
%\item
2. Since $\calF_N^{-1}$ is a transformation which keeps orthonormality, the columns of the matrix calculated as $\bfZ = \calF_N^{-1}\left(\mathrm{orthonorm}(\bfA_g^{-1} \bfR_r)\right)$ also form an orthogonal basis of $\calZ(\bfa)$.
%	\end{enumerate}

\end{remark}

\subsection{Shifting to improve conditioning}
Unfortunately, the circulant matrix $\bfC(\bfa)$ can be rank-deficient; e.g. in the case of the linear series $s_n=c_1 n + c_2$, which
is governed by the GLRR($\bfa$) with $\bfa = (1, -2, 1)^\rmT$. Therefore, instead of solving the linear systems \eqref{eq:lineqMGN},
we consider similar systems with $\bfC(\widetilde{a})$, changing $a$ to $\widetilde{a}$ and then explain how use them to obtain the solutions of \eqref{eq:lineqMGN}.% and \eqref{eq:lineqMGN2}.

Lemma~\ref{lemma:ev_circulant_pre} shows that the eigenvalues of $\bfC(\bfa)$ coincide with the values of the polynomial $g_\bfa(z)$ in nodes of the equidistant grid $\calW = \left\{\exp\big(\frac{\unit 2 \pi j}{N}\big), \; j = 0, \ldots, N-1\right\}$ on the complex unit circle $\spT= \{z \in \spC : |z| = 1 \}$. Therefore, the nondegeneracy of $\bfC(\bfa)$ is equivalent to that there are no roots of the polynomial $g_\bfa(z)$ in $\calW$. The following lemma helps to avoid the problem with zero eigenvalues. Let us define the unitary matrix
\begin{equation}
\label{eq:eqi_grid}
\bfT_M(\alpha) = \diag\left((1, e^{\unit \alpha}, \ldots, e^{\unit (M-1) \alpha})^\rmT\right),
\end{equation}
where $\alpha$ is a real number, $M$ is a natural number.

\begin{lemma}\label{lemma:fourier_1}
	For any real $\alpha$, the following is true: $\bfQ^\rmT(\bfa) \bfx = \bfy$ is attained for some $\bfx \in \spC^N$, $\bfy \in \spC^{N-r}$ if and only if  $\bfQ^\rmT(\tilde \bfa) \left(\bfT_{N}(\alpha)\right) \bfx = \left(\bfT_{N-r}(\alpha)\right) \bfy$, where $\tilde \bfa = \tilde \bfa(\alpha) = \left(\bfT_{r+1}(-\alpha)\right) \bfa$. In addition, the eigenvalues of $\bfC(\tilde \bfa)$ are equal to $g_{\tilde \bfa}(\omega_j) = g_\bfa(\omega_j^{(\alpha)})$, where $\omega_j^{(\alpha)} = \omega_j e^{-\unit \alpha}$.
\end{lemma}
\begin{proof}
	The lemma directly follows from the definitions of the operator $\bfQ(\bfa)$ \eqref{op:Q} and the circulant matrix $\bfC(\bfa)$ \eqref{eq:circulant}.
\end{proof}

The equality $g_{\tilde \bfa}(\omega_j) = g_\bfa(\omega_j^{(\alpha)})$ means that the eigenvalues of $\bfC(\tilde \bfa(\alpha))$ coincide with
the values of the polynomial $g_\bfa(\tilde \omega)$ in $\tilde \omega \in \calW(\alpha) = \left\{\omega_j^{(\alpha)}, \; j = 0, \ldots, N-1\right\}$, where $\omega_j^{(\alpha)} = \exp \left(\unit \left(\frac{2 \pi j}{N} - \alpha \right) \right)$, $\calW(\alpha)$ is the $\alpha$-rotated equidistant grid on $\spT$
(it is sufficient to consider $-\pi/N < \alpha \le \pi/N$, since $\alpha$ and $\alpha +2 \pi/N$ yield the same rotated grid).
Therefore,  $\bfC(\tilde \bfa(\alpha))$ can be made non-degenerate by choosing a suitable  $\alpha$.

\begin{remark} \label{rem:fourier_rotat}
Lemma~\ref{lemma:fourier_1} provides a way for the calculation of an orthonormal basis of $\calZ(\bfa)$. Let us take $\alpha \in \spR$ such that $\bfC(\tilde \bfa)$ is non-degenerate for $\tilde \bfa = \tilde \bfa(\alpha)$.  Using Lemma \ref{lemma:ev_circulant_pre} and Remark~\ref{rem:fourier_orthogon}, we can obtain a matrix $\widetilde \bfZ$ formed from orthonormal basis vectors of $\calZ(\tilde \bfa)$, that is, $\bfQ^\rmT(\tilde \bfa) \widetilde \bfZ = \bm{0}_{(N-r) \times r}$ and $\colspace(\widetilde \bfZ) = \calZ(\widetilde \bfa)$. Then $\bfZ = \left(\bfT_N(-\alpha)\right) \widetilde \bfZ$ has orthonormal columns and $\bfQ^\rmT(\bfa)\bfZ = \bm{0}_{(N-r) \times r}$, that is, $\colspace(\bfZ) = \calZ(\bfa)$.
\end{remark}

In the exact arithmetic, an arbitrary small non-zero value of the smallest eigenvalue of a matrix provides its non-degeneracy. However, in practice, the numerical stability and accuracy of matrix calculations depend on the condition numbers of matrices.
Therefore, the aim of the choice of a proper $\alpha$ is to do the condition number of $\bfC(\tilde \bfa(\alpha))$ as small as possible.
This minimization problem can be approximately reduced to the problem
of maximization of the smallest eigenvalue $|\lambda_\text{min}(\alpha)| = \min_{z \in \calW(\alpha)} | g_\bfa(z) |$ of $\bfC(\tilde \bfa(\alpha))$, since the maximal eigenvalue is not larger than $\max_{z \in \spT} |g_\bfa(z)|$.

\subsection{Algorithm}

By combining Lemmas \ref{lemma:ev_circulant_pre} and \ref{lemma:fourier_1} with Remarks \ref{rem:fourier_orthogon} and \ref{rem:fourier_rotat}, we obtain  Algorithm~\ref{alg:fourier_basis_A} for calculation of an orthonormal basis of $\calZ(\bfa)$.

\begin{algorithm}
	\caption{Calculation of a basis of $\calZ(\bfa)\subset \spC^N$}
	\label{alg:fourier_basis_A}
    \Input{$\bfa \in \spR^r$.}
	\begin{algorithmic}[1]
		\State Find $\alpha_0 = \argmax_{-\pi/N \le \alpha < \pi/N} \min_{z \in \calW(\alpha)} | g_\bfa(z) |$ by means of a 1D numerical optimization method.
		\State Calculate  the vector $\bfa_g = (a_{g, 0}, \ldots, a_{g, N-1})^\rmT$ consisting of the eigenvalues of $\bfC(\widetilde \bfa)$ by $a_{g, j} = g_\bfa\big(\exp(\unit (\frac{2 \pi j}{N} - \alpha_0)\big)$, $j = 0, \ldots, N-1$; $\bfA_g = \diag(\bfa_g)$.
		\State Calculate the matrices $\bfR_r = \calF_N([\bfe_{N-r+1}: \ldots: \bfe_N])$ and  $\bfL_r = \bfA_g^{-1} \bfR_r$.
		\State Find a matrix $\bfU_r \in \spC^{N \times r}$ consisting of orthonormalized columns of the matrix $\bfL_r$
        (e.g, $\bfU_r$ can be obtained by means of the QR decomposition of $\bfL_r$).
		\State Compute $\widetilde \bfZ = \calF_N^{-1}(\bfU_r)$.
		\State \Return $\bfZ = (\bfT_{N}(-\alpha_0)) \widetilde \bfZ \in \spC^{N\times r}$, whose columns form an orthonormal basis of $\calZ(\bfa)$, $\alpha_0$ and $\bfA_g$.
	\end{algorithmic}
\end{algorithm}

%\bigskip
\begin{remark} \label{rem:A_g_C_a}
	Note that the use of the Fourier transform in Algorithm~\ref{alg:fourier_basis_A} allows us to avoid solving the system of linear equations with the matrix $\bfC(\tilde \bfa(\alpha))$. Instead, we invert the diagonal matrix $\bfA_g$, which has the same set of eigenvalues (and, therefore, the same condition number) as the matrix $\bfC(\tilde \bfa(\alpha))$.
\end{remark}

\subsubsection{Numerical properties}
Let us discuss the numerical behavior of Algorithm~\ref{alg:fourier_basis_A}.
The following theorem shows the order of the condition number of the circulant matrix $\bfC(\tilde \bfa(\alpha))$, where $\tilde \bfa(\alpha) = \left(\bfT_{r+1}(-\alpha)\right) \bfa$ is introduced in Lemma~\ref{lemma:fourier_1}, with respect to $\alpha$ in dependence on the series length $N$.
Conventionally, `big O' means an upper bound of the function order, while `big Theta' denotes the exact order.

\begin{theorem} \label{th:gamma}
Let $t$ be the maximal multiplicity of roots of the
polynomial $g_\bfa(z)$ on the unit circle $\spT$. Denote $\lambda_\text{min}(\alpha)$ the minimal eigenvalue of $\bfC (\tilde \bfa(\alpha))$ and $\lambda_\text{max}(\alpha)$ the maximal eigenvalue.
Then
	\begin{enumerate}
		\item for any real sequence $\alpha(N)$, $|\lambda_\text{min}(\alpha)| = O(N^{-t})$;
		\item for any real sequence $\alpha(N)$, $|\lambda_\text{max}(\alpha)| = \Theta(1)$;
		\item there exists such real sequence $\alpha(N)$ that $|\lambda_\text{min}(\alpha)| = \Theta (N^{-t})$.
	\end{enumerate}
\end{theorem}
\begin{proof}
See the proof in Section~\ref{sec:th:gamma}.
\end{proof}

Theorem~\ref{th:gamma} shows that the condition number of the matrix $\bfC(\tilde \bfa(\alpha))$ and $\bfA_g$ used in Algorithm \ref{alg:fourier_basis_A} can be considered as having the order $\Theta(N^{t})$.

\subsubsection{Use of the compensated Horner scheme}
\label{subsec:hornerscheme}
The Horner scheme is an algorithm for evaluating univariate polynomials in floating-point arithmetic. The accuracy of the compensated Horner scheme \cite[Algorithm 4.4 (\code{CompHorner})]{Graillat2008} is similar to the one given by the Horner scheme computed in twice the working precision.

The Horner scheme (we will consider its compensated version) can be directly applied in Algorithm \ref{alg:fourier_basis_A} for calculating the polynomial $g_\bfa$.
Moreover, the Horner scheme can improve the accuracy of the calculation of $\bfU_r$ at step 4 of Algorithm~\ref{alg:fourier_basis_A}; this improvement is important if $\bfL_r$
is ill-conditioned.

To use the advantage of the Horner scheme, let us consider a new way of calculating the matrix $\bfU_r$.
Let $\bfO_r$ be such that $\bfL_r \bfO_r$ consists of orthonormal columns; $\bfO_r$ can be found by either the QR factorization or the SVD.
Then $\bfU_r = \bfA_g^{-1} (\bfR_r \bfO_r)$, where the matrix $\bfR_r$ is calculated at step 3 of Algorithm~\ref{alg:fourier_basis_A}.
Since $(\bfR_r)_{\row{k}} = \left(\exp\big(\frac{\unit 2 \pi r k}{N}\big), \exp\big(\frac{\unit 2 \pi (r-1) k}{N}\big), \ldots, \exp\big(\frac{\unit 2 \pi k}{N}\big)\right)$,
we can reduce the multiplication of $\bfR_r$ by a vector to the calculation of a polynomial of degree $r$ at the point $\exp\big(\frac{\unit 2 \pi k}{N}\big)$.
Therefore, we can accurately calculate the multiplication of $\bfR_r$ by a vector with the help of the Horner scheme. In particular, $\bfR_r \bfO_r$ can be calculated in this way.

\begin{algorithm}
	\caption{Calculation of the basis of $\calZ(\bfa)$ using the Compensated Horner Scheme}
	\label{alga:fourier_basis_A_comp}
    \Input{$\bfa \in \spR^r$.}
	\begin{algorithmic}[1]
		\State Compute $\alpha_0$ and $\bfA_g$ in the same way as at steps 1 and 2 of Algorithm~\ref{alg:fourier_basis_A} except for the use
of the algorithm \code{CompHorner} for calculation of values of the polynomials $g_\bfa$.
		\State Compute $\bfL_r$ and $\bfR_r$ in the same way as at step 3 of Algorithm~\ref{alg:fourier_basis_A}.
		\State Compute $\bfU_r$ in a new manner: find $\bfO_r$ such that $\bfL_r \bfO_r$ consists of orthonormal columns;
calculate $\bfB=\bfR_r \bfO_r$ by means of the algorithm \code{CompHorner}; calculate $\bfU_r = \bfA_g^{-1} \bfB$ directly by matrix multiplication.
\State\Return Matrix $\bfZ$, which is calculated in the same way as at steps 5 and 6 of Algorithm~\ref{alg:fourier_basis_A}, $\alpha_0$, $\bfA_g$.
	\end{algorithmic}
\end{algorithm}
Algorithm~\ref{alga:fourier_basis_A_comp} is a stable analogue of Algorithm~\ref{alg:fourier_basis_A}.

\section{Algorithm for calculation of the projection onto $\calZ(\bfa)$}
\label{sec:proj_our}
Let us describe how to calculate the projection onto $\calZ(\bfa)$ if the basis is given. Together with projection, the algorithm  provides the pseudoinverse to the matrix consisting of the basis vectors.

     We assume that if the matrix $\Sigminus$ is $(2p+1)$-diagonal and positive definite, then it is presented in the form of
  the Cholesky decomposition $\Sigminus = \bfC^\rmT \bfC$; here $\bfC$ is an upper triangular matrix with $p$ nonzero superdiagonals \cite[p. 180]{GoVa13}.
  If $\Sigminus^{-1}$ is $(2p+1)$-diagonal and positive definite, then we consider the representation $\Sigminus = \widehat \bfC^{-1} (\widehat \bfC^{-1})^\rmT$, where $\Sigminus^{-1} = \widehat \bfC^\rmT \widehat \bfC$ is the Cholesky decomposition of $\Sigminus^{-1}$; here $\widehat\bfC$ is an upper triangular matrix with $p$ nonzero superdiagonals.

\begin{remark}\label{rem:wlsinfourier}
	As we mentioned in the beginning of Section~\ref{sec:optim}, the calculation of pseudoinverses ($\inverse{(\bfC \bfZ)}$ or $\inverse{((\widehat \bfC^{-1})^\rmT \bfZ)}$ in our case) can be reduced to solving a linear weighted least-squares problem and therefore their computing can be performed with the help of either the QR factorization or the SVD of the matrix $\bfC \bfZ$ or $(\widehat \bfC^{-1})^\rmT \bfZ$ respectively.
\end{remark}

Algorithm~\ref{alg:proj_calc} provides the algorithm for calculating the pseudo-inverse together with the projection onto a subspace if the basis of this subspace is known.

\begin{algorithm}
	\caption{Calculation of $\winverse{\bfZ}{\bfW}$ and $\Proj_{\bfZ, \bfW} \bfx$ with the use of $\Sigminus = \bfC^\rmT \bfC$ or $\Sigminus^{-1} = \widehat \bfC^\rmT \widehat \bfC$}
	\label{alg:proj_calc}
    \Input{$\bfZ \in \spC^{N\times r}$, $\bfW \in \spR^{N\times N}$ and $\bfx \in \spC^N$.}
	\begin{algorithmic}[1]
		\If{$\Sigminus$ is $(2p+1)$-diagonal}
		\State{Compute the vector $\bfC \bfx$ and the matrix $\bfC \bfZ$.}
		\State{Calculate $\bfq = \inverse{(\bfC \bfZ)} (\bfC \bfx)$, see Remark \ref{rem:wlsinfourier}.}
		\EndIf
		\If{$\Sigminus^{-1}$ is $(2p+1)$-diagonal}
		\State{Compute the vector $(\widehat \bfC^{-1})^\rmT \bfx$ and the matrix $(\widehat \bfC^{-1})^\rmT \bfZ$.}
		\State{Calculate $\bfq = \inverse{((\widehat \bfC^{-1})^\rmT \bfZ)} ((\widehat \bfC^{-1})^\rmT \bfx)$, see Remark \ref{rem:wlsinfourier}.}
		\EndIf
		\State\Return{$\winverse{\bfZ}{\bfW} = \bfq \in \spR^{r\times N}$ and $\Proj_{\bfZ, \bfW} \bfx = \bfZ \bfq  \in \spR^{N}$.}
	\end{algorithmic}
\end{algorithm}

Algorithm~\ref{alg:solution_calc_our} is used for calculating the projection of a given vector.

\begin{algorithm}
	\caption{Calculation of $\Proj_{\calZ(\bfa), \bfW} \tsX$ with the use of special properties of $\calZ(\bfa)$}
	\label{alg:solution_calc_our}
    \Input{$\bfa \in \spR^{r}$, $\bfW \in \spR^{N\times N}$.}
	\begin{algorithmic}[1]
		\State{Compute the matrix $\bfZ(\bfa)$ consisting of basis vectors of $\calZ(\bfa)$ by Algorithm \ref{alg:fourier_basis_A}.}
		\State\Return{Calculate $\Proj_{\bfZ(\bfa), \bfW} \tsX$ by means of Algorithm~\ref{alg:proj_calc}.}
	\end{algorithmic}
\end{algorithm}

\subsection{Computational cost and stability}
\label{sec:comp_cost}
Let us estimate computational costs in flops and study the asymptotic costs as $N \rightarrow \infty$.
The proposed algorithms can be divided into several standard operations with known computational costs. We will use the following asymptotic orders: FFT of a sequence of length $N$ takes $O(N \log N)$ flops \cite[Chapter 1.4.1]{GoVa13}, FFT of a unit vector of length $N$ takes $\Theta(N)$ flops; the Cholesky decomposition of a $(2p+1)$-diagonal matrix $\bfA\in\spR^{N\times N}$ takes $\Theta(N (p+1)^2)$ flops \cite[Chapter 4.3.5]{GoVa13};  solving the  system of linear equations $\bfA\bfx = \bfb$ with $\bfb \in \spR^{N}$ using the obtained decomposition takes additionally $\Theta(N (p+1))$ flops, whereas for $\bfA\bfX = \bfB$ with $\bfB \in \spR^{N\times r}$ the additional cost is $\Theta(N r (p + 1))$ flops; the QR decomposition of an $N\times r$ matrix of rank $r$ takes $\Theta(N r^2)$ flops \cite[Chapter 5.2]{GoVa13}; the pseudo-inversion has the same cost as the QR decomposition, see Remark~\ref{rem:wlsinfourier}; the cost of matrix multiplication is directly determined by their size and structure, in particular, the multiplication of a $(2p+1)$-diagonal $N\times N$ matrix by a vector takes $\Theta(N (p + 1))$ flops \cite[Chapter 1.2.5]{GoVa13}, where $p=0$ corresponds to the case of a diagonal matrix; the computation of a polynomial of order $r$ at $N$ given points takes $\Theta(Nr)$ flops.

\paragraph{Algorithm~\ref{alg:fourier_basis_A}}
Although the implementations of Algorithm~\ref{alg:fourier_basis_A} differ for the case when $\Sigminus^{-1}$ is $(2p+1)$-diagonal
and the case when $\Sigminus$ is $(2p+1)$-diagonal, the asymptotic computational cost is the same. Algorithm~\ref{alg:fourier_basis_A} includes computing the $N \times N$ diagonal matrix $\bfA_g$, where each diagonal value is obtained using the calculation of a polynomial of order $r$ (step 2); solving a system of linear equations given by a diagonal matrix (step 3); FFT of $r$ unit vectors (step 3); FFT of $r$ arbitrary vectors (step 5); the QR decomposition (step 4); the multiplication of a diagonal matrix by a vector $r$ times (step 6). The search of optimal rotations at step 1 of Algorithm~\ref{alg:fourier_basis_A} serves for increasing of the algorithm stability. Therefore we can fix the number of iterations in this search. Since the computational cost of calculating the objective function is  $O(N r)$ flops, the cost of step 1 is also $O(N r)$ flops. Therefore, Algorithm~\ref{alg:fourier_basis_A} requires $O(r N \log N + N r^2)$ flops, or $O(N \log N)$ for a fixed $r$.

\paragraph{Algorithm~\ref{alg:proj_calc}}
Calculating the projection by Algorithm~\ref{alg:proj_calc} includes the multiplication by a $p$-diagonal matrix and the QR decomposition for the pseudoinverse computation that leads to $\Theta(N r^2 + N rp)$ operations.

\paragraph{Algorithm~\ref{alg:solution_calc_our}}
Algorithm~\ref{alg:solution_calc_our} consists of the calls of Algorithms~\ref{alg:fourier_basis_A} and~\ref{alg:proj_calc}. 
Therefore, for either $\Sigminus$ or $\Sigminus^{-1}$ is $(2p+1)$-diagonal, the asymptotical computational cost of Algorithm~\ref{alg:solution_calc_our} is $O(r N \log N + N r^2+Npr)$ flops, or $O(N \log N+ Np)$ for a fixed $r$.

\paragraph{Stability}
The main ``stability bottlenecks'' of Algorithms~1--4 is the inversion of the matrix $\bfA_g$ in Algorithm~\ref{alg:fourier_basis_A}.
The inversion of the matrix $\bfA_g$ serves for solving the linear systems \eqref{eq:lineqMGN}
% and \eqref{eq:lineqMGN2}
 in a stable and fast way (see Remark~\ref{rem:A_g_C_a}).
 Thus, let us discuss the  orders of the condition numbers of the matrix $\bfA_g$ as the time-series length $N$ tends to infinity.

 Theorem~\ref{th:gamma} shows that the order of the condition number of the matrix $\bfA_g$ is $\Theta(N^{t})$, where $t$ is the maximal multiplicity of roots of the characteristic polynomial  $g(\bfa)$ \eqref{eq:pol_z} on the unit circle.

%\clearpage
\section{Calculation of projections and the VPGN algorithm}
\label{sec:MGNand VPGN}

The algorithm VPGN described in Section~\ref{sec:VPGN} contains calculating the projection $\Proj_{\calZ(\bfa), \bfW} \bfx$. Let us discuss how the method of projecting suggested in Section~\ref{sec:proj_our} influences the VPGN implementation in comparison with the implementation suggested in \cite{Usevich2014}.

\subsection{Calculating the projections in \cite{Usevich2014}}
Let us describe the algorithm described in Section~\ref{sec:VPGN} in the form suggested in \cite{Usevich2014}.
 In \cite{Usevich2014}, calculating the projection $\Proj_{\calZ(\bfa), \bfW} \bfx$ is performed by means of constructing the projection $\Proj_{\calQ(\bfa), \bfW}$ onto the orthogonal compliment $\calQ(\bfa)$ and then subtracting from the identity matrix: $\bfI_N - \Proj_{\calQ(\bfa), \bfW}$.
Thus, the following relation is used in \cite{Usevich2014}:
\begin{equation} \label{eq:spZa_kostya}
	\Proj_{\calZ(\bfa), \bfW} \bfx = \left( \bfI_N - \bfW^{-1} \bfQ(\bfa) \bm\Gamma^{-1}(\bfa) \bfQ^\rmT(\bfa) \right) \bfx,
	\end{equation}
	where $\bm\Gamma(\bfa) = \bfQ^\rmT(\bfa) \Sigminus^{-1} \bfQ(\bfa)\in \spR^{(N-r)\times (N-r)}$ (see Lemma~\ref{th:varproj}).
The calculation of $\Proj_{\calZ(\bfa), \bfW}$ by \eqref{eq:spZa_kostya} needs computing the matrix  $\bm \Gamma^{-1}(\bfa)$.
Below we write down Algorithm~\ref{alg:solution_calc_vp}, which was used in the paper \cite{Usevich2014}, with a fast computation of $\bm \Gamma(\bfa)$ and its inverse (see Algorithm~\ref{alg:gamma_inverse}). Algorithm~\ref{alg:gamma_inverse} uses the matrix $\widehat \bfC$, which is defined as at the beginning of Section~\ref{sec:proj_our}, i.e. $\Sigminus^{-1} = \widehat \bfC^\rmT \widehat \bfC$ is the Cholesky decomposition of $\Sigminus^{-1}$.

\begin{algorithm}
	\caption{Calculation of $\bm\Gamma^{-1}(\bfa) \bfv$ by the method from \cite{Usevich2014}}
	\label{alg:gamma_inverse}
    \Input{$\bfa \in \spR^{r}$, $\bfW \in \spR^{N\times N}$, $\bfW^{-1}$ is $(2p+1)$-diagonal ($p\leq N$), $\bfv \in \spR^{N-r}$}
	\begin{algorithmic}[1]
		\State{Calculate the matrix $\widehat \bfC \bfQ(\bfa)$, which has $m+1$ non-zero diagonals, where $m=\min(p+r, N-r)$.}
		\State{Calculate $(2m+1)$-diagonal matrix $\bm\Gamma(\bfa) = (\widehat \bfC \bfQ(\bfa))^\rmT(\widehat \bfC \bfQ(\bfa))$.}
		\State{Calculate the Cholesky decomposition $\bm\Gamma(\bfa) = (\bm\Gamma_\mathrm{c})^\rmT \bm\Gamma_\mathrm{c}$, where $\bm\Gamma_\mathrm{c}$ is $(m+1)$-diagonal.}
		\State\Return $\bm\Gamma^{-1}(\bfa)\bfv = \bm\Gamma_\mathrm{c}^{-1}\left( (\bm\Gamma_\mathrm{c}^\rmT)^{-1}\bfv \right)$
	\end{algorithmic}
\end{algorithm}

Algorithm~\ref{alg:gamma_inverse} is used for calculating the projection $\Proj_{\calZ(\bfa), \bfW}$ in Algorithm~\ref{alg:solution_calc_vp} in the way similar to that in \cite{Usevich2014}.

\begin{algorithm}
	\caption{Calculation of $\Proj_{\calZ(\bfa), \bfW} \tsX$ by the method from \cite{Usevich2014} using \eqref{eq:spZa_kostya}}
	\label{alg:solution_calc_vp}
	\Input{$\bfa \in \spR^{r}$, $\bfW \in \spR^{N\times N}$.}
	\begin{algorithmic}[1]
		\State{Compute $\bfv = \bfQ^\rmT(\bfa) \tsX \in \spR^{N-r}$}
		\State{Compute $\bfy = \bm\Gamma^{-1}(\bfa) \bfv \in \spR^{N-r}$ using Algorithm~\ref{alg:gamma_inverse}}
		\State\Return{$\Proj_{\calZ(\bfa), \bfW} \tsX = \tsX - \bfW^{-1}\bfQ^\rmT(\bfa)\bfy$}
	\end{algorithmic}
\end{algorithm}

\subsubsection{Computational cost and stability}
\paragraph{Algorithm~\ref{alg:gamma_inverse}}
Let $\Sigminus^{-1}$ be $(2p+1)$-diagonal.
Algorithm~\ref{alg:gamma_inverse} includes computing the Cholesky factorization of a ($2m+1$)-diagonal matrix of order $(N-r) \times (N-r)$, where $m \le p + r$ (step 3); solving a system of linear equations using the obtained decomposition (step 4); the multiplications of matrices with $p+1$ and $r+1$ non-zero diagonals (step 1), $(m+1)$ and $(m+1)$ non-zero diagonals (step 2). This gives us the asymptotic cost $O(N r^2 + N p^2)$ flops or $O(N p^2)$ for a fixed $r$.

For the case when $\Sigminus$ is $(2p+1)$-diagonal, $p > 0$, there is no implementation of Algorithm~\ref{alg:gamma_inverse} faster than with cubic (in $N$) asymptotic complexity, since $\bm \Gamma(\bfa)$ (see Section~\ref{sec:MUdetails}) is not a banded matrix.

\paragraph{Algorithms~\ref{alg:solution_calc_vp}}
Besides the call of Algorithm~\ref{alg:gamma_inverse}, Algorithm~\ref{alg:solution_calc_vp} includes the multiplications by a $(r+1)$-diagonal matrix (step 1, step 3)
and a $(2p+1)$-diagonal matrix (step 3). Therefore, the asymptotic cost is the same as for Algorithm~\ref{alg:gamma_inverse}.

\paragraph{Stability}
The main ``stability bottlenecks'' of Algorithm~\ref{alg:solution_calc_vp} is solving the systems of linear equations with matrices related to $\bfQ(\bfa)$.
For Algorithm~\ref{alg:solution_calc_vp}, it is the matrix $\bm\Gamma(\bfa)$ whose inversion is constructed in Algorithm~\ref{alg:gamma_inverse}. Let us discuss the  order of the condition number of this matrix as the time-series length $N$ tends to infinity.

We consider the case when $\bfW^{-1}$ is banded, since otherwise the computational cost of Algorithm~\ref{alg:solution_calc_vp} is very large.
For fast inversion, the diagonals of the matrix $\bm\Gamma(\bfa)$ are computed explicitly; then the Cholesky factorization is used. It is shown in \cite[Section 6.2]{Usevich2014} that the condition number of $\bm\Gamma(\bfa)$ is $O(N^{2t})$. (Compare with the condition number $O(N^{t})$ of the ``stability bottlenecks'' matrix in the proposed Algorithm~\ref{alg:solution_calc_our}, which is discussed in Section~\ref{sec:comp_cost}.)

Certainly, the inversion of $\bm\Gamma(\bfa)$ can be performed with better stability.
For example, one can use the QR factorization of the matrix $\bfW^{-1/2} \bfQ(\bfa)$ instead of the inversion of $\bm\Gamma(\bfa)$. However, the QR factorization does not exploit the banded structure of matrix $\bm\Gamma(\bfa)$, therefore, it is significantly slower than the Cholesky factorization if $\bfW^{-1}$ is banded.

%\clearpage
\subsection{Algorithms VPGN and S-VPGN}
\label{sec:VPGN_alg}
Algorithm~\ref{alg:gauss_newton_vp} contains the formal description of the VPGN and S-VPGN algorithms with different implementations of the projection $\Proj_{\calZ(\bfa), \bfW} \bfx$.

Recall notation: $\fullop$ and $S_\tau^\star(\Ai0) = \Proj_{\calZ(
\bfa), \bfW} \tsX$, where $\bfa = \fullop(\Ai0)$, are introduced in Sections~\ref{sec:param} and \ref{sec:VPGN} respectively.

\begin{algorithm}
	\caption{Variable Projection Gauss-Newton method (the VPGN and S-VPGN versions)}
	\label{alg:gauss_newton_vp}
    \Input{$\tsX \in \spR^N$, $\bfa_0\in \spR^{r+1}$, a stopping criterion STOP.}
	\begin{algorithmic}[1]
		\State{Set $k = 0$,  $\bfb^{(0)} = \bfa_0$.}
		\Repeat{}
		\State{Choose $\i0$ such that $b_\i0^{(k)}\neq 0$; for example, find $\i0 = \argmax_{i} |b_i^{(k)}|$. Calculate $\bfa^{(k)} = c\bfb^{(k)}$, where $c$ is such that $c b_\tau^{(k)} = -1$, and take $\Ai0^{(k)} = \fullop^{-1}(\bfa^{(k)})$.}
		\State{Calculate $\tsS_k  = S_\tau^\star(\Ai0^{(k)})$  using Algorithm~\ref{alg:solution_calc_vp} (VPGN) or Algorithm~\ref{alg:solution_calc_our} (S-VPGN) to compute $\Proj_{\calZ(\bfa^{(k)}) , \bfW}$.}
		\State{Calculate $\bfJ_{S_\tau^\star}(\Ai0^{(k)})$ by \eqref{eq:vpformula} applying Algorithm~\ref{alg:gamma_inverse} to compute $\bm\Gamma^{-1}(\bfa^{(k)}) \bfv$ for $\bfv = \bfQ^\rmT(\bfe_j) \Proj_{\calZ(\bfa^{(k)}), \bfW} \tsX$ and for  $\bfv = \bfQ^\rmT(\bfa^{(k)}) \tsX$.}
		\State{Calculate
			$\Delta_k = \winverse{\bfJ_{S_\tau^\star}(\Ai0^{(k)})}{\bfW} (\tsX - \tsS_k)$} by applying Algorithm~\ref{alg:proj_calc} to computing the pseudoinverse.
		\State{Perform a step of size $\gamma_k$ for the line search in the descent direction given by $\Delta_k$ using Algorithm~\ref{alg:solution_calc_vp} (VPGN) or Algorithm~\ref{alg:solution_calc_our} (S-VPGN) to compute $\Proj_{\calZ(\bfa^{(k)}) , \bfW}$ to calculate $S_\tau^\star$. For example, find $0\le\gamma_k\le 1$ such that
          $\|\tsX - S_{\tau}^\star(\Ai0^{(k)} + \gamma \Delta_k) \|_{\bfW} \le \|\tsX - S_\tau^\star(\Ai0^{(k)}) \|_{\bfW}$
          by the backtracking method \cite[Section 3.1]{nocedal2006numerical}.}
		\State{Set $\Ai0^{(k+1)} = \Ai0^{(k)} + \gamma_k \Delta_k$, $\bfb^{(k+1)} = \fullop(\Ai0^{(k+1)})$.}
		\State{Set $k = k+1$.}
		\Until{STOP}
		\State\Return{$\widetilde \tsS = S_{\i0}^\star(\Ai0^{(k)})$ as an estimate of the signal.}
	\end{algorithmic}
\end{algorithm}

In both algorithms, the calls of Algorithms~\ref{alg:solution_calc_vp} and Algorithm~\ref{alg:solution_calc_our} are supplemented by the calls of Algorithm~\ref{alg:gamma_inverse} at step 4.
Therefore, if $\Sigminus^{-1}$ be $(2p+1)$-diagonal, S-VPGN has asymptotic computational cost $O(r N \log N + N r^2 + Nrp + N p^2) = O(r N \log N + N r^2  + N p^2)$, instead of $O(N r^2 + Np^2)$ for VPGN, that is, S-VPGN is slightly slower.
For the case when $\Sigminus$ is $(2p+1)$-diagonal, $p > 0$, the complexity of both algorithms is $O(N^3)$.

Thus, we consider the case when $\Sigminus^{-1}$ be $(2p+1)$-diagonal to compare algorithms by stability. We expect that the use of the more stable Algorithm~\ref{alg:solution_calc_our} instead of Algorithm~\ref{alg:solution_calc_vp} improves the stability of Algorithm~\ref{alg:gauss_newton_vp}.

\begin{remark}
 Algorithm~\ref{alg:solution_calc_our} calls Algorithm~\ref{alg:fourier_basis_A}. If the compensation Horner scheme is used in Algorithm~\ref{alg:fourier_basis_A} (see Section~\ref{subsec:hornerscheme}), we call the optimization algorithm S-VPGN-H.
\end{remark}

\subsection{Numerical comparison of stability} \label{subsec:speed}
\subsubsection{Example}
\label{subsec:basisacc}
With the help of Lemma~\ref{lemma:locminnec}, we construct an example, where a local solution of \eqref{eq:wls} is known.
For constructing a solution of rank $r=3$, we use the well-known theory about the relation of linear recurrence relations,
characteristic polynomials, their roots and the explicit form of the series, see e.g. the book \cite[Sections 3.2]{Golyandina2013} with a brief description of this relation in the context of time series structure.

Let $\tsY_N^\star = (b t_1^2, \ldots, b t_N^2)^\rmT$, where $t_i$, $i=1,\ldots,N$, form the equidistant grid in $[-1; 1]$ and the constant $b$ is such that $\|\tsY_N^\star\|=1$. The series $\tsY_N^\star$ satisfies the GLRR($\bfa^*$) for $\bfa^* = (1, -3, 3, -1)^\rmT$. Since the last component of $\bfa^*$ is equal to $-1$, we can say that the series satisfies the LRR($\bfa^*$).
Denote $\widehat \tsR_N = (c |t_1|, \ldots, c |t_N|)^\rmT$, where the constant $c$ is such that $\|\widehat \tsR_N\| = 1$.
Construct the observed series as $\tsX_N = \tsY_N^\star + \tsR_N$, where $\tsR_N = \widehat \tsR_N - \Proj_{\calZ((\bfa^*)^2), \Sigminus} \widehat \tsR_N$.
Thus, the pair $\tsX_0 = \tsY_N^\star$ and $\tsX = \tsX_N$ satisfies the conditions of Lemma \ref{lemma:locminnec}, which provides the necessary conditions for local minima. The sufficient condition (the positive definiteness of the Hessian matrix of the objective function  $\|\tsX -  \tsS(\Si0, \Ai0)\|^2_\bfW$ \cite[Theorem 2.3]{nocedal2006numerical}) was tested numerically for $N < 100$.

\paragraph{Implementation}
In practice, we should generate the time series $\tsX_N=(x_1,\ldots,x_N)^\rmT$ with high numerical precision which is enough for comparing
the algorithms, which solve the problem \eqref{eq:wls}, by their accuracy.
The main difficulty lies in calculating the projection $\Proj_{\calZ(\bfa_0^2), \Sigminus}$.
The GLRR($\bfa_0$) with $\bfa_0 = (1, -3, 3, -1)^\rmT$ corresponds to the characteristic polynomial $g_{\bfa_0}(t)=(t-1)^3$ with the coefficients taken from $\bfa_0$.
Therefore, the GLRR($\bfa_0^2$) corresponds to the characteristic polynomial
$g^2_{\bfa_0}(t)=(t-1)^6$ and a basis of $\calZ(\bfa_0^2)$ consists of polynomials of degree not greater than $5$.
To obtain
the projection, we use the Legendre polynomials \cite{Belousov2014} of degree from $0$ to $5$, which are calculated at the points $t_i$ as a basis of $\calZ(\bfa_0^2)$.
Then the constructed basis is orthogonalized.

\subsubsection{Comparison of projection accuracy}

Before comparing the solution stability, let us compare the accuracy of different methods implementing the projection onto the subspaces $\calZ(\bfa)$.
The accuracy of the projection calculation is important for solving the problem \eqref{eq:wls}, since the constructed solution should belong to $\calD_r$ consisting of subspaces $\calZ(\bfa)$ for different $\bfa$.

Consider the time series $\tsP_N = \tsY^{*}_N + \tilde \tsR_N$, where $\tilde \tsR_N = \widehat \tsR_N - \Proj_{\calZ(\bfa^*), \Sigminus} \widehat \tsR_N$ and calculate $\Proj_{\calZ(\bfa^*), \Sigminus} \tsP_N$.
The comparison is performed for projection implementations done by the methods VP (Algorithm~\ref{alg:solution_calc_vp}), S-VP (Algorithm~\ref{alg:solution_calc_our}) and S-VP-H (Algorithm~\ref{alg:solution_calc_our} with the use of the compensated Horner scheme) for different $N$ from $20$ to $50000$.

As a measure of accuracy, we take the Euclidean distance from  $\Proj_{\calZ(\bfa^*), \Sigminus} \tsP_N$ to $\tsY^{*}_N$, which theoretically equals zero. For the correctness of numerical results, the values of the constructed time series $\tsP_N$ were calculated with the best available accuracy; since the subspace $\calZ(\bfa^*)$ consists of polynomials of degree not larger than 2, the basis of $\calZ(\bfa^*)$ is calculated with the help of Legendre polynomials, in the same way as in the implementation of the example.

\begin{figure}[!hbt]
	\centering
	\includegraphics[width=7cm]{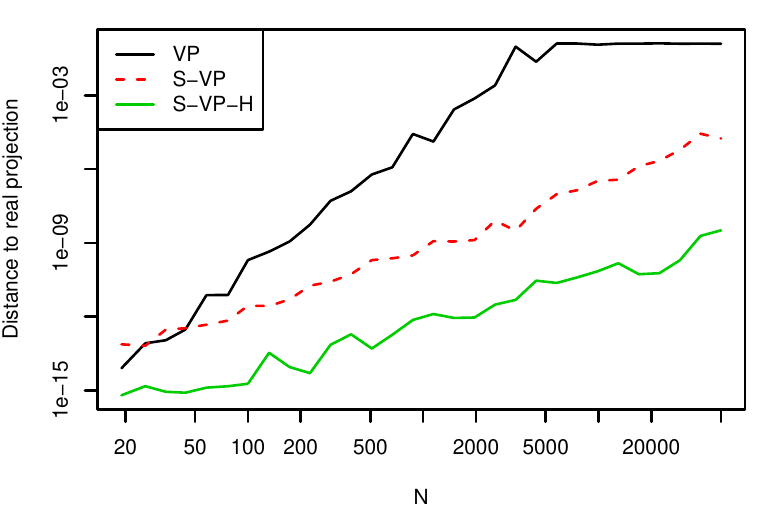}
	\caption{Comparison of projection accuracy, for different $N$.}
	\label{fig:fig_proj}
\end{figure}

For simplicity, consider the non-weighted case, when $\Sigminus$ is the identity matrix.
The results are presented in Fig.~\ref{fig:fig_proj}. On can see that the S-VP and S-VP-H methods have much smaller errors than the VP method, which fails for N larger 2000. The compensated Horner scheme considerably improves the accuracy of projecting.

\subsubsection{Comparison of solution stability}
\paragraph{Line search and stopping criteria}
The compared algorithms contain a line search in the descent direction $\Delta_k$.
 The line search method and the stopping criteria are not specified in the algorithms.
 Let us provide details concerning the implementation of the line search at step 7 and the stopping criterion in Algorithms \ref{alg:gauss_newton_vp}.
We implemented the backtracking line search method \cite[Section 3.1]{nocedal2006numerical} in the direction $\Delta_k$ starting from the step size $\gamma = 1$ (the full step) and then dividing $\gamma$ by 2. The backtracking stops when
\begin{equation}
\label{eq:searchstop}
\|\tsX - \tsS^\star(\Ai0^{(k)} + \gamma \Delta_k) \|_{\bfW} \le \|\tsX - \tsS^\star(\Ai0^{(k)}) \|_{\bfW};
\end{equation}
 then $\gamma_k = \gamma$. If there is no such $\gamma$ for $\gamma = 1, 1/2, 1/4, \ldots, 2^{-50}$, then we set $\gamma_k = 0$.
  The stopping criterion of the whole algorithm is the equality $\gamma_k = 0$, which means that the current iteration can not improve the approximation to the solution.

\begin{figure}[!hbt]
	\centering
	(a)\includegraphics[width=7cm]{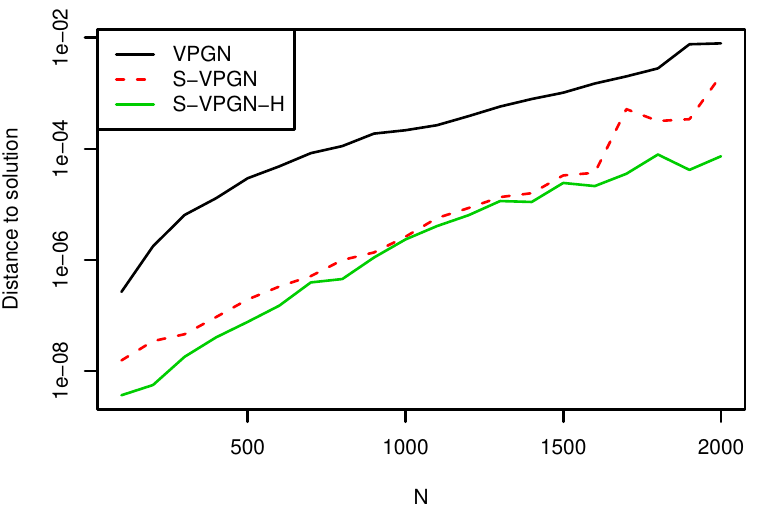}
	(b)\includegraphics[width=7cm]{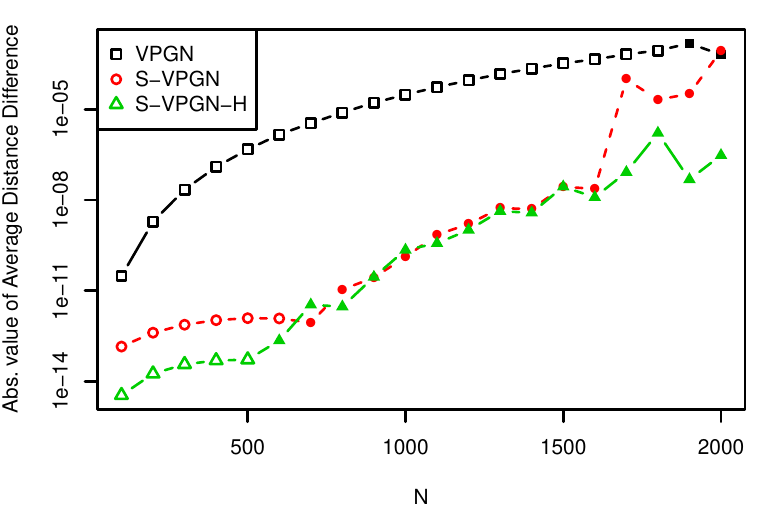}
	\caption{Comparison of algorithms (a) by distance to the solution and (b)  by absolute differences between the values of the objective function; for different $N$.}
	\label{fig:kostya_comp_disp}
\end{figure}

\paragraph{Comparison}
 Denote $\widetilde \tsY^\star$ the result of an algorithm participating in the comparison.
The main comparison was done by accuracy, that is, by the Euclidean distance between $\widetilde \tsY^\star$ and the solution $\tsY^\star_N$ (Fig.~\ref{fig:kostya_comp_disp}(a)).

In addition, the algorithms were compared by discrepancy between the values of the objective function at the final
point of the algorithm and at the point of local minimum, i.e. by $\|\tsX_N -  \widetilde \tsY^\star\| - \| \tsX_N - \tsY^\star_N\|$ (Fig.~\ref{fig:kostya_comp_disp}(b)).

The algorithms were started from the GLRR($\bfa_0$), where  $\bfa_0 =\bfa^* + 10^{-6} \bfd^\rmT$ and
each components of $\bfd$ is randomly distributed in $[-1,1]$.
We used 100 simulations to obtain the average results.

Figure~\ref{fig:kostya_comp_disp} shows that the accuracy of S-VPGN and S-VPGN-H is better than that of VPGN, for window lengths  $N$ less than approximately several thousand when the algorithms start to fail, although the S-VP and S-VP-H projections are still working (see Fig.~\ref{fig:fig_proj}). The possible reason is that step 5 of Algorithm~\ref{alg:gauss_newton_vp} was not improved. The advantage of S-VPGN-H over S-VPGN is visible but weak.

\section{Conclusion}
\label{sec:conclusion}
The properties of the subspace of low-rank time series were studied; these properties are based on the chosen family of local parameterizations related to generalized linear recurrent relations GLRR($\bfa$). It was shown that this parameterization is smooth and therefore allows one to consider different numerical optimization methods (e.g. the Gauss-Newton method) for solving least-squares problems.
We proved (Theorem~\ref{th:tangent}) that the tangent subspace at the point $\tsS$, which is governed by a GLRR($\bfa$), can be described in terms of the GLRR($\bfa^2$).
This fact allows one to construct first-order linear approximations to functions at points from $\calD_r$.

Then, in Section~\ref{sec:ZofA} we present a numerically stable algorithm of projecting a series onto the set $\calZ(\bfa)$ of time series, which are governed by the GLRR($\bfa$). The computational cost of the proposed projection algorithm is $O(N \log N + Np)$, if either the weight matrix $\bfW$ or its inverse is $(2p+1)$-diagonal. This algorithm can be useful for numerical solutions of different approximation problems related to the SLRA problems; we demonstrate it by means of improving the stability of the known algorithm of low-rank time series approximation from \cite{Usevich2014}.

%\clearpage
%\setcounter{page}{1}
\appendix
%\renewcommand{\thesection}{A}
%\section{Appendix}
\label{sec:app}

\section{Additional information}
\subsection{Correspondence between notations}
\label{sec:notation}
For the convenience of comparisons, in Table~\ref{table:defines} we present the correspondence between the notation used in this paper and the notation from  \cite{Usevich2014, Usevich2012}.

\begin{table}[!hh] 	\centering
	\caption{Correspondence between notations}
    \label{table:defines}
	\begin{tabular}{|c|c|c|c|c|c|c|c|c|}
		\hline
		This paper & $\tsX$ & $N$ & $r+1$ & $1$ & $N-r$ & $\bfa$ & $\Sigminus$ & $\bm\Gamma(\bfa)$ \\ \hline
		Usevich \& Markovsky & $p_D$ & $n_p$ & $m$ & $d$ & $n$ & $R$ & $W$ & $\Gamma$\\
		\hline
	\end{tabular}
\end{table}

\subsection{Rank of \eqref{eq:model}}
\label{sec:rank_calc}
\begin{proposition}
Let a series $\tsS$ of length $N$ have the form \eqref{eq:model}, $0\le \omega_k \le 0.5$, $0\le \phi_k < 2\pi$ and $m_k\ge 0$ be the polynomial degree, $k=1,\ldots,d$.
Suppose that the pairs $(\alpha_k, \omega_k)$ are different.
Also, assume that if $\omega_k = 0$ or $\omega_k = 0.5$, then $\phi_k \neq 0$.
Let $r_k$ be equal to 2 if $0 < \omega_k < 0.5$ and be equal to 1 otherwise.
Then the rank of $\tsS$ is equal to $\sum_{k=1}^d (m_k+1)r_k$ for sufficiently large $N$.
\end{proposition}

\begin{proof}
The assertion about the ranks of real-valued time series is the consequence of the analogous results for complex-valued time series.
If a series $\tsC=(c_1,\ldots,c_N)$ has terms $c_n = \sum_{k=1}^s P_{m_k}(n) \mu_k^n$ with different complex $\mu_k$, then its rank $r$ is
equal to  $\sum_{k=1}^s (m_k+1)$. This directly follows from the explicit form of the basis of the column space of the trajectory matrix $\calT_{r+1}(\tsC)$, which consists of  $r$ linearly independent vectors
$(1^i \mu_k^1, 2^i \mu_k^2, \ldots (r+1)^i \mu_k^{r+1})^\rmT$, $k=1,\ldots,s$, $i=0,\ldots,m_k$.
The rank of a real-valued time series is induced by the presentation of $\exp(\alpha_k n) \sin(2\pi \omega_k n + \phi_k)$, $0<\omega_k < 0.5$, as a linear combination
of $\mu^n$ and $\overline{\mu}^n$, where $\mu = \exp(\alpha_k + \unit 2\pi \omega_k)$ and $\overline{\mu}$ is the complex conjugate to $\mu$.
\end{proof}

\subsection{Lemma about $\overline{\calD_r}$}
\label{sec:closure}
\begin{lemma}
	$\tsS \in \overline{\calD_r}$ if and only if $\tsS$ is governed by a GLRR($\bfa$) defined by a vector $\bfa \in \spR^{d+1}$, $d \le r$.
\end{lemma}
\begin{proof}
Let us consider the set of matrices $\calM_{\le r} \subset \spR^{L \times K}$ of rank not larger than $r$, and $\calM_{=r}$ the set of matrices of rank  $r$. Fix $L = r+1$, $K = N - L + 1 = N-r$.
Denote $\widehat \calD_r = \{\tsX\in \spR^N : \rank \calT_{r+1}(\tsX) \le r \} = \calT_{r+1}^{-1}(\calM_{\le r} \cap \calH) = \bigcup_{s=1}^r \calD_{s}$.
By definition,  $\calD_r = \{\tsX : \rank \calT_{r+1}(\tsX) = r \} = \calT_{r+1}^{-1}(\calM_{=r} \cap \calH)$.

It is known that $\overline{\calM_{=r}} = \calM_{\le r}$, see \cite{Lewis2008}. Thus, we have $\overline{\calD_r} = \overline{\calT_{r+1}^{-1}(\calM_{=r} \cap \calH)} = \calT_{r+1}^{-1}(\overline{\calM_{=r} \cap \calH}) \subset \calT_{r+1}^{-1}(\overline{\calM_{=r}} \cap \calH) = \calT_{r+1}^{-1}(\calM_{\le r} \cap \calH) = \widehat \calD_r$.

	To prove $\overline{\calD_r} = \widehat \calD_r$, we show that any $\tsS \in \widehat \calD_r$ can be approximated by a series $\tsX \in \calD_r$ with arbitrary precision. Let $\tilde{r} < r$ and $\tsS \in \calD_{\tilde{r}}$ satisfy a GLRR($\tilde{\bfa}$), $\tilde{\bfa} = (a_1, \ldots, a_{\tilde{r}+1})^\rmT \in \spR^{\tilde{r} + 1}$. It is sufficient to show that we can approximate $\tsS$ by $\tsX \in \calD_{\tilde{r} + 1}$; then we can obtain an approximating series from $\calD_r$ by subsequent approximations with ranks increased by 1.

Let us take such real $\mu$ that the series $\tsD=(\mu,\mu^2,\ldots,\mu^N)^\rmT$ of rank 1 is not governed by the GLRR($\tilde{\bfa}$). Denote $\bfd_M = (\mu,\mu^2,\ldots,\mu^M)^\rmT$; then $\tsD=\bfd_N$. For any real $\alpha \ne 0$, we have
	$\tsX(\alpha) = \tsS + \alpha \tsD \in \widehat \calD_{\tilde{r} + 1}$, since the series $\tsX(\alpha)$ is governed by the GLRR($\bfb$) with $\bfb = (\mu a_1, \mu a_2 - a_1, \mu a_3 - a_2, \ldots, \mu a_{\tilde{r}+1} - a_{\tilde{r}}, -a_{\tilde{r} + 1})^\rmT \in \spR^{\tilde{r} + 2}$. Thus, $\rank \tsX(\alpha) \le \tilde{r} +1$.
	
	Now let us show that $\rank \tsX(\alpha) \ge \tilde{r} +1$.
 We need to show that $\rank \calT_{\tilde{r} + 1}(\tsX(\alpha)) = \tilde{r}+1$ for any $\alpha \ne 0$.
 Due to \cite[Corollary 8.1]{Marsaglia1974}, it is enough to show that the column and row spaces of $\calT_{\tilde{r} +1}(\tsS)$ and $\calT_{\tilde{r} +1}(\alpha \tsD)$ have empty intersection.
 We know that $\colspace\left(\calT_{\tilde{r} +1}(\alpha \tsD)\right) = \sspan(\bfd_{\tilde{r} +1})$ and $\rowspace\left(\calT_{\tilde{r} +1}(\alpha \tsD)\right) = \sspan(\bfd_{N-\tilde{r}})$. Also, note that a vector $\bfv$ belongs to $\colspace\left({\calT_{\tilde{r} + 1}(\tsS)}\right)$ if and only if $\tilde{\bfa}^\rmT \bfv = 0$, and a vector $\bfu$ belongs to $\rowspace\left({\calT_{\tilde{r} + 1}(\tsS)}\right)$ if and only if
  $\left(\bfQ^{N-\tilde{r}, \tilde{r}}(\tilde{\bfa})\right)^\rmT \bfu= \bm{0}_{N-2\tilde{r}}$. However, by construction of $\tsD$, $\tilde{\bfa}^\rmT \bfd_{\tilde{r} +1} \ne 0$, and $\left(\bfQ^{N-\tilde{r}, \tilde{r}}(\tilde{\bfa})\right)^\rmT \bfd_{N-\tilde{r}}  \ne \bm{0}_{N-2\tilde{r}} $. Therefore, we have
 $\rowspace\left(\calT_{\tilde{r} + 1}(\tsS)\right) \cap \rowspace\left(\calT_{\tilde{r} + 1}(\alpha \tsD)\right) = \emptyset$ and $\colspace\left( \calT_{\tilde{r} + 1}(\tsS)\right) \cap \colspace\left(\calT_{\tilde{r} + 1}(\alpha \tsD)\right) = \emptyset$. The lemma is proved, since $\alpha$ can be an arbitrarily small positive number.
\end{proof}

    \section{Proofs of propositions from the paper}
    \subsection{Proof of Theorem~\ref{th:parameterization} and Proposition~\ref{prop:parameterization}}
    \label{sec:th:parameterization}

    \begin{proof}
    The first statement of Proposition~\ref{prop:parameterization} will provide the parameterizing mapping introduced in Theorem~\ref{th:parameterization}
    if we prove the correctness of \eqref{eq:param} and \eqref{eq:param_rev}, the uniqueness of $S_\tau$ satisfying relations of Theorem \ref{th:parameterization}, then prove that $S_\tau$ is an injective mapping and \eqref{eq:param_rev} defines the inverse of the mapping $S_\tau$ given in \eqref{eq:param}.

    Let us prove the correctness of \eqref{eq:param}.
    To begin with, we show that the matrix $\bfZ_{\row{\calI({\i0})}}$ is not singular and therefore invertible. This will be a consequence of non-singularity of $(\bfZ_0)_{\row{\calI({\i0})}}$ for any basis of $\calZ(\bfa_0)$.

    Let us represent $\bfa_0$ as $\bfa_0 = (0, \ldots, 0, b_{r_m +1 }, \ldots, b_{1}, 0, \ldots, 0)^\rmT\in \spR^{r+1}$, with $r_b$ zeroes at the beginning and $r_e$ zeroes at the end, $r_e + r_b + r_m = r$. Let us construct a matrix $\bfZ_0^\star= [\bfZ_\text{begin}: \bfZ_\text{middle}:\overline{} \bfZ_\text{end}]$ consisting of three blocks: $\bfZ_\text{begin} = \begin{pmatrix}
    \bfI_{r_b} \\
    \bm{0}_{(N - r_b)\times r_b}
    \end{pmatrix}$, $\bfZ_\text{middle} = \begin{pmatrix}
    \bm{0}_{r_b \times r_m} \\
    \widehat \bfZ_\text{middle} \\
    \bm{0}_{r_e \times r_m}
    \end{pmatrix} $, $\bfZ_\text{end} = \begin{pmatrix}
    \bm{0}_{(N - r_e)\times r_e} \\
    \bfI_{r_e}
    \end{pmatrix}$, where the columns of the matrix $\widehat \bfZ_\text{middle} \in \spR^{(N - r_b - r_e) \times r_m}$ form a basis of the space of time series of length $N - r_b - r_e$ governed by the LRR with coefficients $-b_2/b_1, \ldots, -b_{r_m +1 }/b_1$. Since $\{1, \ldots, r_b\} \cup \{N - {r_e}+1, \ldots, N\} \subset \calI(\i0)$ and any submatrix of size $r_m \times r_m$ of $\widehat \bfZ_\text{middle}$ is non-degenerate \cite[Prop. 2.3]{Usevich2010}, we obtain the non-degeneracy of $(\bfZ_0^\star)_{\row{\calI({\i0})}}$. Any other matrix which consists of basis vectors of $\calZ(\bfa_0)$ can be represented in the form $\bfZ_0^\star \bfP$ with a non-singular matrix $\bfP \in \spR^{r\times r}$. Therefore, matrix $(\bfZ_0^\star \bfP)_{\row{\calI({\i0})}}$ is also non-degenerate.

    Now let us prove the non-degeneracy of $\bfZ_{\row{\calI({\i0})}}$. Since $\calZ(\bfa)$ is the orthogonal complement to $\calQ(\bfa)$, $\Proj_{\calZ(\bfa)}$ can be represented as a continuous function $\Proj_{\calZ (\bfa)} = \bfI_N - \Proj_{\bfQ(\bfa)}$ of $\bfa$, $\bfa \ne \bm{0}_{r+1}$, where $\bfQ(\bfa)$ is defined in \eqref{op:Q}.
   Note that the determinant of $\bfZ_{\row{\calI({\i0})}}$ is a continuous function of $\bfZ$. In turn, $\bfZ$ continuously depends on $\Ai0$.
   Since $\bfZ(\bfa_0) = \bfZ_0$ and the determinant of $(\bfZ_0)_{\row{\calI({\i0})}}$ is non-zero, there is a neighborhood of $(\bfa_0)_{\calK(\i0)}$, such that
   the determinant of $\bfZ_{\row{\calI({\i0})}}$ is not zero; therefore, the matrix $\bfZ_{\row{\calI({\i0})}}$ is invertible.

    The constructed mapping \eqref{eq:param} does not depend on $\bfZ_0$. Indeed, for any non-singular matrix $\bfP\in \spR^{r\times r}$:
    $\left(\Proj_{\calZ (\bfa)} \bfZ_0 \bfP\right) \left((\Proj_{\calZ (\bfa)} \bfZ_0 \bfP)_{\row{\calI({\i0})}}\right)^{-1} = \bfZ \left(\bfZ_{\row{\calI({\i0})}}\right)^{-1}$.

	Let us demonstrate that the properties of $S_\tau$, which are stated in Theorem~\ref{th:parameterization}, are fulfilled; i.e., show that $\tsS \in \calD_r$, the series $\tsS$ satisfies the GLRR($\bfa$) and $(\tsS)_{\vecrow{\calI(\i0)}} = \Si0$. The series $\tsS$ satisfies the GLRR($\bfa$), since each column of the matrix $\bfZ$ satisfies the GLRR($\bfa$).
	To prove that $\tsS \in \calD_r$, consider the matrix $\calT_{r+1}(\tsS_0)$ and choose a submatrix of size $r \times r$ with
non-zero determinant. Then take the submatrix $\bfB$ of the matrix $\calT_{r+1}(S_\tau(\Si0, \Ai0))$ with the same location. Its determinant is
 a continuous function of $(\Si0, \Ai0)$, since the function given in \eqref{eq:param} is continuous. Therefore, there exists a neighborhood
 of $\left((\tsS_0)_{\calI(\i0)}, (\bfa_0)_{\calK(\i0)}\right)^\rmT$, where the determinant of $\bfB$ is non-zero; thus, $\tsS \in \calD_r$.
The condition $(\tsS)_{\vecrow{\calI(\i0)}} = \Si0$ is fulfilled, since
    \begin{equation*}
    (\tsS)_{\vecrow{\calI(\i0)}} =
    \left(\bfZ_{\row{\calI({\i0})}} \left(\bfZ_{\row{\calI({\i0})}}\right)^{-1}\right) \Si0 = \Si0.
    \end{equation*}

    Let us explain the uniqueness of the mapping $S_\tau$ satisfying the relations of Theorem~\ref{th:parameterization}.
    Let  $\widehat S$ be a different mapping satisfying the relations of Theorem~\ref{th:parameterization}, $\widehat \tsS = \widehat S(\Si0, \Ai0) \in \calD_r$. We know that $\widehat \tsS \in \calZ(\bfa)$. Therefore, columns of $\bfZ$ contain a basis of $\calZ(\bfa)$. Let $\widehat \tsS = \bfZ \bfv$ and $\bfv \in \spR^r$ be the coefficients of the expansion of $\widehat \tsS$ in the columns of $\bfZ$. Then the following is fulfilled: $(\bfZ \bfv)_{\calI({\i0})} = \Si0$. However, $\bfZ_{\row{\calI({\i0})}} \bfv = \Si0$ together with the invertibility of $\bfZ_{\row{\calI({\i0})}}$ leads to $\bfv =  \left(\bfZ_{\row{\calI({\i0})}}\right)^{-1} \Si0$. Therefore, $\widehat \tsS = \left(\bfZ \left(\bfZ_{\row{\calI({\i0})}}\right)^{-1}\right) \Si0 = \tsS$.

    Let us prove that $S_\tau$ is an injective mapping. We choose two different sets of parameters $\big(\Si0^{(1)}, \Ai0^{(1)}\big)^\rmT$, $\big(\Si0^{(2)}, \Ai0^{(2)}\big)^\rmT$ in the vicinity of $\left((\tsS_0)_{\calI(\i0)}, (\bfa_0)_{\calK(\i0)}\right)^\rmT$ and consider $\tsX_1 = S_\tau\big(\Si0^{(1)}, \Ai0^{(1)}\big)$, $\tsX_2 = S_\tau\big(\Si0^{(2)}, \Ai0^{(2)}\big)$. If $\Si0^{(1)} \ne \Si0^{(2)}$, then $\tsX_1 \ne \tsX_2$, since $(\tsX_1)_{\vecrow{\calI(\i0)}} \ne (\tsX_2)_{\vecrow{\calI(\i0)}}$. Let $\Si0^{(1)} = \Si0^{(2)}$ be fulfilled, but $\Ai0^{(1)} \ne \Ai0^{(2)}$. This means that the orthogonal complements $\sspan(\fullop(\Ai0^{(1)}))$ and $\sspan(\fullop(\Ai0^{(2)}))$ to $\colspace{\left(\calT_{r+1}(\tsX_1)\right)}$ and $\colspace{\left(\calT_{r+1}(\tsX_2)\right)}$ respectively are different and therefore these column spaces differs.
    Thus, $\tsX_1\ne \tsX_2$.

    Let us prove the correctness of \eqref{eq:param_rev}.
    According to the statement of Proposition~\ref{prop:parameterization}, $\Ai0$ defined in \eqref{eq:param_rev} is obtained from a renormalization of $\hat \bfa = \hat \bfa(\tsS)$ such that the $\i0$-th element becomes equal to $-1$.
    Let us prove the correctness of this definition of $\Ai0$, i.e., the possibility to renormalize $\hat \bfa$. Consider the matrix $\bfS = \calT_{r+1}(\tsS) \in \spR^{(r+1)\times (N-r)}$. Let $\calJ$ be a subset of indices such that the submatrix $(\bfS_0)_{\col{\calJ}}\in \spR^{(r+1)\times r}$ has rank $r$, where $\bfS_0 = \calT_{r+1}(\tsS_0)$. Then $\Proj_{\calL(\tsS)}$ can be represented as a continuous function $\Proj_{\calL(\tsS)} = \Proj_{\bfS_{\col{\calJ}}}$ in the vicinity of $\tsS_0$; therefore, we can choose a neighborhood of $\tsS_0$ in which $\hat a_\i0$ does not vanish.

    Let us explain that $\eqref{eq:param_rev}$ gives the inverse of the mapping $S_\tau$. Let $\tsS = S_\tau(\Si0, \Ai0)$. The values $\Si0 = (\tsS)_{\vecrow{\calI({\i0})}}$ are taken directly from the time series. The series $\tsS$ is governed by the GLRR($\hat \bfa$) since the vector $\hat \bfa$ is orthogonal to $\colspace(\calT_{r+1}(\tsS))$ by its definition. But the series $\tsS$ is governed by the GLRR($\bfa$); hence, $\bfa$ coincides with $\hat \bfa$ up to normalization. Therefore, renormalization of $\hat \bfa$ gives us the required $\Ai0$. This consideration concludes the proof.
    \end{proof}

\subsection{Proof of Theorem~\ref{th:gamma}}
\label{sec:th:gamma}
\begin{proof}
Denote by $\upangle{x}{y}$ the angle between two points on the complex unit circle $\spT$, $0 \le \upangle{x}{y} \le \pi$. Let us prove the first statement. Consider a root $z_1 \in \spT$ of multiplicity $t_1$, $t_1\le t$, of the polynomial $g_\bfa(z)$; then for any $\alpha$ we have $\min_{w \in \calW(\alpha)} \upangle{w}{z_1} \le \frac{\pi}{N}$ by the Dirichlet principle. Let us fix any $0 \le \alpha_0 < 2 \pi$ and choose $w_0 = \argmin_{w \in \calW(\alpha_0)} \upangle{w}{z_1}$. Since $|z_1 - w_0|=O(1/N)$, we have $|\lambda_\text{min}(\alpha)| \le |g_\bfa(w_0)| = O(N^{-t})$.
	
	To prove the second statement, let us find any point $x \in \spT$ for which $|g_\bfa(x)| = \max_{z \in \spT} |g_\bfa(z)| > 0$ is fulfilled. Again, by the Dirichlet principle, we have $\min_{w \in \calW(\alpha)} \upangle{w}{x} \le \frac{\pi}{N}$ for any $\alpha$. Let us choose $w_1 = \argmin_{w \in \calW(\alpha_0)} \upangle{w}{x}$. Since $|x - w_1|=O(1/N)$ and $g_\bfa(z)$ is continuous, we have $|\lambda_\text{max}(\alpha)| \ge |g_\bfa(w_1)| = \Omega(1)$, which with $|\lambda_\text{max}(\alpha)| = O(1)$ proves the second part.
	
	To prove the third statement, let us construct a piecewise approximation of $g_\bfa(z)$ in $z$. Consider the decomposition $g_\bfa(z) = p_\bfa(z) q_\bfa(z)$, where the roots of $p_\bfa(z)$ belong to $\spT$ while the roots of $q_\bfa(z)$ do not. By construction, $\inf_{z \in \spT} |q_\bfa(z)| > 0$.
	
	Let $z_1, \ldots, z_k$ be the roots of $p_\bfa(z)$ with multiplicities $t_1, \ldots, t_k$. We split the circle $\spT$ into $k$ semi-open non-intersecting arcs $\calS_1, \ldots, \calS_k$, $\spT = \bigcup_{1 \le i \le k} \calS_i$, such that $z_i \in \calS_i$ for any $i$ and $z_j \notin \overline{\calS_i}$ for any $j \ne i$ ($\overline{\calS_i}$ denotes the closure of $\calS_i$), which leads to $\inf_{z \in \calS_i} \left| {p_\bfa(z)}/{(z-z_i)^{t_i}} \right| > 0 $.
	
	To finish the proof, we need to show that there exists $0 \le \alpha = \alpha(N) < 2 \pi$
	such that $$\min\limits_{w \in \calW(\alpha(N)), \; 1 \le i \le k} \upangle{w}{z_i} = \Theta(1/N).$$
	Denote for $0 \le \mu < \pi / N$ and $z \in \spT$
	$$\calB_{z, \mu} = \{0 \le \alpha < 2 \pi : \min_{w \in \calW(\alpha)} \upangle{w}{z} \le \mu \}.$$
	The set $\calB_{z, \mu}$ has the explicit form:
	$$\calB_{z, \mu} = \bigcup_{0 \le j \le N-1} \Big\{ \Arg \left( \exp \left( \unit \left( \frac{2 \pi j}{N} + y \right)  / z \right) \Big|_{-\mu \le y \le \mu} \right) \Big\}.$$

	Let us comment this expression. Consider $\omega_j^{(\alpha)} = \exp \left(\unit \left(\frac{2 \pi j}{N} - \alpha \right) \right)$ and choose $\alpha_j$ such that $\upangle{\omega_j^{(\alpha_j)}}{z} \le \mu$. This means that the polar angle of the ratio $\omega_j^{(\alpha_j)}/z$ belongs to the interval $[-\mu, \mu]$, i.e. $\omega_j^{(\alpha_j)}/z \in \{ \exp \left( \unit x \right) |_{-\mu \le x \le \mu} \}$. Evidently,
$$\exp \left(\unit \left({2 \pi j}/{N} - \alpha_j \right) \right) \in \{ z \exp \left( \unit x \right) |_{-\mu \le x \le \mu} \},$$
is equivalent to
 $$\exp \left(\unit \left(\alpha_j - {2 \pi j}/{N} \right) \right) \in \{ \exp \left( \unit y \right) / z |_{-\mu \le y \le \mu} \},$$
  where $y = -x$. Finally, note that $\alpha_j \in  \calM_j = \left\{ \Arg (\exp \left( \unit ({2 \pi j}/{N} + y) \right) / z ) |_{-\mu \le y \le \mu} \right\}$. The inequality $\min_{w \in \calW(\alpha)} \upangle{w}{z} \le \mu$ is valid if $\alpha$ is equal to one of $\alpha_0, \ldots, \alpha_j$.
Therefore, the union of all such sets $\calM_j$ for $j = 0, \ldots, N-1$ gives us $\calB_{z, \mu}$.
	
	The Lebesgue measure of $\calB_{z, \mu}$ is equal to $\mes \calB_{z, \mu} = 2 \mu N$ for $\mu < \pi/N$.
	Let us take $\mu = \frac{\pi}{2Nk}$ and consider $\calB = \bigcup_{1 \le i \le k} \calB_{\mu, z_i}$.
	Since $\mes \calB \le \pi$, we obtain $\mes \widehat \calB \ge \pi$ for $\widehat \calB = [0; 2 \pi) \setminus \calB$,
	which means that $\widehat \calB$ is not the empty set.
	Thus, we have proved that for any $\alpha \in \widehat \calB$
	$$\min_{w \in \calW(\alpha), \; 1 \le i \le k} \upangle{w}{z_i} > \frac{\pi}{2Nk}.$$
	
	Let us fix an arbitrary $\alpha_0 \in \widehat \calB$ and consider any $w \in \calW(\alpha_0)$. %Let $i$ be such that $w \in \calS_i$.
	For each $i$ such that $w \in \calS_i$, $|w - z_i| = \Theta(1/N)$.
	Then $|g_\bfa(w)| = |q_\bfa(w)| \left| \displaystyle{\frac{p_\bfa(w)}{(w-z_i)^{t_i}}} \right| |(w-z_i)^{t_i}| \ge C \Theta(N^{-t_i})$,
	where $C > 0$ is some constant.
\end{proof}

\end{document}